\documentclass[11pt]{amsart}
\usepackage{texdraw}
\usepackage{amsmath, amssymb}

\newcommand{\vs}{\vskip 12 pt}
\newcommand{\ra}{\rightarrow}

\newcommand{\PP}{\mathbb P}
\newcommand{\BC}{\mathbb C}
\newcommand{\BZ}{\mathbb Z}

\newcommand{\CB}{\mathcal B}
\newcommand{\CR}{\mathcal R}

\newcommand{\CF}{\mathcal F}
\newcommand{\CS}{\mathcal S}

\newcommand{\BS}{\mathbb S}
\newcommand{\ba}{\mathbf a}
\newcommand{\bb}{\mathbf b}
\newcommand{\bc}{\mathbf c}
\newcommand{\bI}{\mathbf I}

\newcommand{\fm}{\mathfrak{m}}
\newcommand{\fn}{\mathfrak{n}}
\newcommand{\fp}{\mathfrak{p}}

\newcommand{\fl}{\mathfrak{l}}

\newtheorem{prop}{Proposition}[section]
\newtheorem{thm}[prop]{Theorem}
\newtheorem{lem}[prop]{Lemma}

\theoremstyle{remark}
  \newtheorem{rk}[prop]{Remark}
\theoremstyle{definition}
 \newtheorem{exam}[prop]{Example}
 \newtheorem{defn}[prop]{Definition}

\newenvironment{pf}[1]{\noindent {\it{Proof.}} {#1}}{$ \square $ \vs}

\newsavebox{\bx}
\savebox{\bx} (60,40)[bl]{
 \put(0,0){\line(1,0){10}}
 \put(10,0){\line(0,1){10}}
 \put(10,10){\line(-1,0){10}}
 \put(0,10){\line(0,-1){10}} }

\begin{document}

\title{Rigidity of singular Schubert varieties in $Gr(m,n)$}

\author{Jaehyun Hong}

\begin{abstract}
Let $\ba=(p_1^{q_1}, \cdots, p_r^{q_r})$ be a partition and
$\ba'=({p_1'}^{q_1'}, \cdots, {p_r'}^{q_r'})$ be its conjugate.
 We will
prove that if $q_i, q_i \geq 2$ for all $1 \leq i\leq r$, then
any irreducible subvariety $X$ of $Gr(m,n)$ whose homology class
is an integral multiple of the Schubert class $[\sigma_{\ba}]$ of
type $\ba$ is a Schubert variety of type $\ba$.
\end{abstract}

\maketitle


\section{Introduction}

Let $Gr(m,n)$ be the Grassmannian of $m$-planes in $\BC^n$. For a
partition $\ba=(a_1, \cdots, a_m)$ a Schubert variety
$\sigma_{\ba}$ of type $\ba$ is defined by the set of all
$m$-planes $E$ such that $ \dim\, (E \cap \BC^{n-m+i-a_i} )\geq i,
1 \leq i \leq m$ for a choice of a flag $\{ \BC^1 \subset \cdots
\subset \BC^n \}$. Then they form a basis for the homology space
$H_*(Gr(m,n), \BZ)$.

For $\ba=(p^q)^*=((n-m)^{m-q}, (n-m-p)^q)$, the Schubert variety
$\sigma_{\ba}$ of type $\ba$ is smooth and they are all the smooth
Schubert varieties in $Gr(m,n)$. The {\it Schur rigidity} of
smooth Schubert varieties in $Gr(m,n)$ is proved in \cite{W},
\cite{B} and \cite{Ho}: for any smooth Schubert variety
$\sigma_{\ba}$ in $Gr(m,n)$ other than a non-maximal linear space,
any irreducible subvariety whose homology class is an integral
multiple of the Schubert class $[\sigma_{\ba}]$ of type $\ba$ is a
Schubert variety of type $\ba$.

In this paper we will prove the Schur rigidity of singular
Schubert varieties of certain types in $Gr(m,n)$.

\vs {\bf Theorem.} {\it Let $\ba=(p_1^{q_1}, \cdots, p_r^{q_r})$
be a partition and let $\ba'=({p'_1}^{q'_1}, \cdots,
{p'_s}^{q'_r})$ be its conjugate.  Then the Schubert variety
$\sigma_{\ba}$ is Schur rigid if
 $q_i, q'_i \geq 2$ for all $i \leq r$.
}

\vs

The proof is divided by two parts  as in \cite{W}, \cite{B} and
\cite{Ho}: Schubert rigidity and the equality $\CB_{\ba}
=\CR_{\ba^*}$. Define the Schubert differential system $\CB_{\ba}$
by the differential system with a fiber at $x \in Gr(m,n)$ given
by the set of all the tangent space of Schubert varieties of type
$\ba$ passing through $x$. If any irreducible integral variety of
$\CB_{\ba}$ is a Schubert variety of type $\ba$, then we say that
the Schubert variety $\sigma_{\ba}$ is {\it Schubert rigid}.

Putting together  the tangent space of all the subvarieties $X$
with $[X]=r[\sigma_{\ba}]$, $r \in \BZ$ at each point,  we get
another differential system $\CR_{\ba^*}$, which we call the Schur
differential system(\cite{B}, \cite{W}). By the construction, if
any irreducible integral variety of $\CR_{\ba^*}$ is a Schubert
variety of type $\ba$, then  the Schubert variety $\sigma_{\ba}$
is  Schur rigid. Furthermore, the Schubert differential system
$\CB_{\ba}$ is always contained in the Schur differential system
$\CR_{\ba^*}$. Thus the equality $\CB_{\ba}=\CR_{\ba^*}$ and the
Schubert rigidity is a necessary and sufficient condition for the
Schur rigidity(\cite{B}). While proving the equality
$\CB_{\ba}=\CR_{\ba^*}$ is computing integral elements of exterior
differential systems, which is an algebraic problem, proving the
Schubert rigidity is finding integral varieties of a differential
system, which is a local differential geometric problem.

\vs
 When the Schubert variety $\sigma_{\ba}$ is singular, we will
express the Schubert variety $\sigma_{\ba}$ as the locus of a
family of smooth Schubert varieties $\sigma_{\bb}$ of type $\bb$
on $Gr(m,n)$ parameterized by a Schubert variety $\sigma_{\bc}$ of
type $\bc$ in another Grassmannian which is the parameter space of
the Schubert variety of type $\bb$ on $Gr(m,n)$. Then the Schubert
variety $\sigma_{\ba}$ is Schubert rigid if the Schubert
varietyies $\sigma_{\bb}$ and $\sigma_{\bc}$ are Schubert rigid.
So we can use the Schubert rigidity of smooth Schubert varieties
and the induction to prove the Schubert rigidity of singular
Schubert varieties in $Gr(m,n)$ of  types in  Theorem(Section 3).

To prove the equality $B_{\ba}=R_{\ba^*} \subset Gr(k, E^* \otimes
Q)$ for the fibers of $\CB_{\ba}$ and $\CR_{\ba^*}$ at $[E] \in
Gr(m,n)$, we use the description $R_{\ba}=Gr(|\ba|, E^* \otimes Q)
\cap \PP(\BS_{\ba}(E^*) \otimes \BS_{\ba'}(Q))$ as in \cite{Ho}
and then  we compute the complement of the tangent space of
$B_{\ba}$ in the tangent space of $Gr(k, E^* \otimes Q)$ by hands,
while in \cite{Ho} the theory of Lie algebra cohomology developed
by Kostant used to compute it(Section 4).

\vs {\bf Acknowledgements} I wish to thank Professor D. Burns for
suggesting the idea of using a natural foliation by smooth
Schubert varieties on singular Schubert varieties to prove their
rigidity. Some results of the present paper have been worked out
while I was visiting MSRI in Berkeley in November 2003.
I would like to thank this institute for its hospitality, and
Professor R. Bryant, one of the organizers of the program
``Differential Geometry'', for stimulating discussions and
encouragement. This work was supported by the Post-Doctorial
Fellowship Program of Korea Science and Engineering
Foundation(KOSEF).


\section{Differential systems}

\subsection{Schubert differential systems}

Let $Gr(m,n)$ be the Grassmannian of $m$-dimensional subspaces of
$V=\BC^n$.
  Let $P(m,n)$ be the set $$\{ \ba=(a_1, \cdots, a_m) |\,\, n-m
\geq a_1 \geq \cdots \geq a_m \geq 0 \}$$ of partitions. Fix a
flag $\{ V_{\bullet} \}$ of $V$ with $\dim V_i = i$. For $\ba \in
P(m,n)$, define the {\it Schubert variety}
$\sigma_{\ba}(V_{\bullet})$ {\it of type} $\ba$ by the set $$\{ E
\in Gr(m,n)|\,\, \dim(E \cap V_{n-m+i-a_i}) \geq i \}.$$ Then
$\sigma_{\ba}(V_{\bullet})$ is a subvariety of $Gr(m,n)$ of
codimension $|\ba|:=a_1 +\cdots + a_m$. By varying the flag $\{
V_{\bullet}\}$, we get a family of Schubert varieties
$\sigma_{\ba}$ of type $\ba$.

Let $\ba \in P(m,n)$. Define its {\it dual} $\ba^*$ by
$$\ba^*=(n-m-a_m, \cdots, n-m-a_1)$$
\noindent and  define its {\it conjugate} $\ba'=(a'_1, \cdots,
a'_{n-m})$  by $$a'_i=\sharp \{ j | a_j \geq i \} \text{ for } 1
\leq i \leq n-m.$$  \noindent
 The {\it Young diagram} $Y_{\ba}$ is defined by the set of
 boxes consisting of $a_i$ boxes in the $i$-th row, the row of
 boxes lined up on the left.
 \noindent Then the Young diagram $Y_{\ba^*}$ is obtained by rotating
 the complement of the Young diagram $Y_{\ba}$  by $180$ degree and
 the Young diagram $Y_{\ba'}$  is obtained by transposing
  the Young diagram $Y_{\ba}$.

\begin{exam}$m=5$, $n-m=6$. Young diagrams $Y_{\ba}, Y_{\ba^*}, Y_{\ba'}$
for $\ba=(6,6,4,2,2)$.  $\ba^*=(4,4,2)$, $\ba'=(5,5,3,3,2,2)$

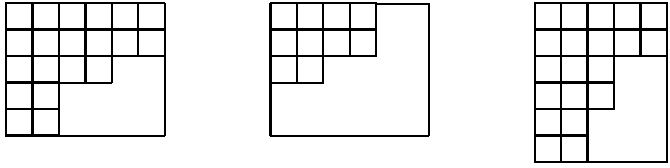
\vs \vs

\begin{picture}(300, 80)

\put(0,60){\usebox{\bx}} \put (10,60){\usebox{\bx}}
\put(20,60){\usebox{\bx}} \put (30,60){\usebox{\bx}}
\put(40,60){\usebox{\bx}} \put (50,60){\usebox{\bx}}

\put(0,50){\usebox{\bx}} \put (10,50){\usebox{\bx}}
\put(20,50){\usebox{\bx}} \put (30,50){\usebox{\bx}}
\put(40,50){\usebox{\bx}} \put (50,50){\usebox{\bx}}

\put(0,40){\usebox{\bx}} \put (10,40){\usebox{\bx}}
\put(20,40){\usebox{\bx}} \put (30,40){\usebox{\bx}}

\put(0,30){\usebox{\bx}} \put (10,30){\usebox{\bx}}

\put(0,20){\usebox{\bx}} \put (10,20){\usebox{\bx}}

\put(20,20){\line(1,0){40}} \put(60,20){\line(0,1){30}}

\put(100,60){\usebox{\bx}} \put (110,60){\usebox{\bx}}
\put(120,60){\usebox{\bx}} \put (130,60){\usebox{\bx}}

\put(100,50){\usebox{\bx}} \put (110,50){\usebox{\bx}}
\put(120,50){\usebox{\bx}} \put (130,50){\usebox{\bx}}

\put(100,40){\usebox{\bx}} \put (110,40){\usebox{\bx}}

\put(100,20){\line(1,0){60}} \put(160,20){\line(0,1){50}}
\put(100,20){\line(0,1){20}} \put(160,70){\line(-1,0){20}}
\put(200,60){\usebox{\bx}} \put (210,60){\usebox{\bx}}
\put(220,60){\usebox{\bx}} \put (230,60){\usebox{\bx}}
\put(240,60){\usebox{\bx}}

\put(200,50){\usebox{\bx}} \put (210,50){\usebox{\bx}}
\put(220,50){\usebox{\bx}} \put (230,50){\usebox{\bx}}
\put(240,50){\usebox{\bx}}

\put(200,40){\usebox{\bx}} \put (210,40){\usebox{\bx}}
\put(220,40){\usebox{\bx}}

\put(200,30){\usebox{\bx}} \put (210,30){\usebox{\bx}}
\put(220,30){\usebox{\bx}}

\put(200,20){\usebox{\bx}} \put (210,20){\usebox{\bx}}

\put(200,10){\usebox{\bx}} \put (210,10){\usebox{\bx}}

\put(220,10){\line(1,0){30}} \put(250,10){\line(0,1){50}}
\end{picture}

\vs

\end{exam}

Define $\fn_{\ba}$ to be the vector space of matrices $Z=(z^a_i)
\in M_{n-m, m}$ that satisfy $z^a_i=0$ when $a >n-m-a_i$. This is
the tangent space of the Schubert variety  $\sigma_{\ba}$ and is
obtained by rotating
 the Young diagram $Y_{\ba^*}$ by 90 degree in counterclockwise.

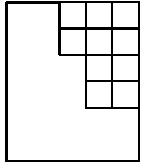
\vs

\begin{picture}(300, 80)

\put(100,70){\line(-1,0){20}} \put(80,70){\line(0,-1){60}}
\put(80,10){\line(1,0){50}} \put(130,10){\line(0,1){60}}

\put(100,60){\usebox{\bx}} \put (110,60){\usebox{\bx}}
\put(120,60){\usebox{\bx}}

\put(100,50){\usebox{\bx}} \put (110,50){\usebox{\bx}}
\put(120,50){\usebox{\bx}}

\put(110,40){\usebox{\bx}} \put (120,40){\usebox{\bx}}

\put(110,30){\usebox{\bx}} \put (120,30){\usebox{\bx}}

\end{picture}

\vs

We will use the notation $\ba=(p_1^{q_1}, \cdots, p_r^{q_r})$,
$p_r \not=0$
 for the partition with $q_1$ $p_1$'s, $\cdots$, $q_r$ $p_r$'s.

\begin{defn} Let $M$ be a manifold and let $Gr(k, TM)$ be the
Grassmannian bundle of $k$-subspaces of the tangent bundle $TM$. A
subvariety $\CF$ of $Gr(k, TM)$ is called a {\it differential
system} on $M$. A subvariety $X$ of $M$ is said to be an {\it
integral variety} of the differential system $\CF$ if at each
smooth point $x \in X$, the tangent space $T_xX$ is an element of
the fiber $\CF_x$. We say that $\CF$ is {\it integrable} if at
each point $x \in M$ and $y \in \CF_x$, there is an integral
variety passing through $x$ and tangent to the subspace $W_y$ of
$T_xM$ corresponding to $y$.
\end{defn}

\begin{defn} For each $\ba \in P(m,n)$, the {\it Schubert differential
system} $\CB_{\ba}$ of type $\ba$ is the differential system with
a fiber consisting of the tangent space to the Schubert varieties
$\sigma_{\ba}$ of type $\ba$ passing through a given point.  We
say that $\sigma_{\ba}$ is {\it Schubert rigid} if Schubert
varieties of type $\ba$ are the only irreducible integral
varieties of $\CB_{\ba}$. If, furthermore, there is a unique
Schubert variety passing through a given point and tangent to a
given tangent subspace, then we say that $\CB_{\ba}$ is {\it
strongly rigid}.
\end{defn}

\begin{rk}
For a smooth Schubert variety in $Gr(m,n)$, the Schubert
differential system is strongly rigid ([W], Proposition
\ref{smooth}) and this will play a central role in proving the
Schubert rigidity of  singular Schubert varieties.
\end{rk}

Let $[E] \in Gr(m,n)$ and $Q=\BC^n /E$. Then $T_{[E]}Gr(m,n)=E^*
\otimes Q$ and $\wedge^k(E^* \otimes Q)=\oplus_{|\ba|=k}
\BS_{\ba}(E^*) \otimes \BS_{\ba'}(Q)$, where $\BS_{\ba}$ is the
Schur functor of type $\ba$. There exists a $SU(n)$-invariant
positive $(k,k)$-form $\phi_{\ba}$ which can be written as the sum
$(\sqrt{-1})^{k^2}\sum_i \xi_i \wedge \overline{\xi}_i$ at $[E]$,
where $\{\xi_i \}$ is an orthonormal basis of $\BS_{\ba}(E)
\otimes \BS_{\ba'}(Q^*) \subset \wedge^k (T^*_{[E]}Gr(m,n))$. Then
we have $\int_{\sigma_{\ba^*}} \phi_{\bb}=\delta_{\ba}^{\bb}$ for
$\ba, \bb \in P(m,n)$(For details, see \cite{B}).

\begin{defn} For each $\ba \in P(m,n)$, the
{\it Schur differential system} $\CR_{\ba}$ of type $\ba$ is
defined by the intersection
$$\cap_{\bb \not=\ba, |\bb|=|\ba|} Z(\phi_{\bb}),$$
\noindent where $Z(\phi_{\bb})$ is the set of $|b|$-subspace of
$T_E(Gr(m,n))$ on which $\phi_{\bb}$ vanishes. We say that
$\sigma_{\ba}$ is {\it Schur rigid} if Schubert varieties  of type
$\ba^*$ are the only integral varieties of $\CR_{\ba}$. \end{defn}

\begin{rk} \label{rk2}

(1) $\CB_{\ba^*}$ is contained in $\CR_{\ba}$, and $X $ is a
subvariety of $Gr(m,n)$ with $[X]=r[\sigma_{\ba^*}]$ for an
integer $r$ if and only if $X$ is an integral variety of
$\CR_{\ba}$(\cite{B}, \cite{W}).

(2) $\CB_{\ba^*}$ is closed and $\CR_{\ba}$ is connected (Remark 2
and Remark 12 of \cite{B}).
\end{rk}


\subsection{$F$-structures and integral varieties}
\label{Fstructure}

\begin{defn} Let $F $  be a submanifold of $ Gr(k, V)$
with a transitive action of a subgroup of $GL(V)$. A fiber bundle
$\CF \subset Gr(k, TM)$ on a manifold $M$ of dimension $n=\dim V$
is said to be an {\it $F$-structure}   if  at each point $x \in M$
there is a linear isomorphism $\varphi(x): V \ra T_xM$ such that
the induced map $\varphi(x)^k : Gr(k, V) \ra Gr(k, T_xM)$ sends
$F$ to $\CF_x$.

\end{defn}

Schubert differential systems are $F$-structures for various
$F$'s.  Integrability and the uniqueness can be obtained by
studying the cohomology space $H^{k,1}(F)$ and $H^{k,2}(F)$
associated to $(W_f, V, T_fF)$ for $f \in F$(\cite{G}). When an
$F$-structure is integrable, if $H^{k,1}(F)=0$, then the family of
all integral varieties passing through a fixed point and tangent
to a fixed subspace has dimension $\sum_{j \leq k-1}\dim
H^{j,1}(F)$. In particular, if $H^{1,1}(F)=0$, there is only one
such integral subvariety. Higher cohomology gives the information
on the higher jet of the integral varieties.

In general, Schubert differential system for singular Schubert
variety has order $\geq 2$, i.e. its integral varieties are
determined by higher jets. But smooth Schubert varieties depend
only 1-jet and there is a canonical map from the differential
system of a singular Schubert variety to that of a certain smooth
Schubert variety(Proposition \ref{surj}). So, in this paper, we
will consider only the first cohomology $H^{1,1}(F)$, which
contains the information on the 2-jets of integral varieties.


\begin{defn} Let $F $  be a subvariety of $ Gr(k, V)$
with a transitive action of a subgroup of $GL(V)$. Let  $\CF$ be
an $F$-structure on $M$ with the projection map $\pi:\CF \ra M$.
For $x \in M$ and $y \in \CF_x$, let $W_y$ denote the $k$-subspace
of $T_xM$ corresponding to $y$. For a $k$-subspace $H$ of $T_y\CF$
such that $\pi_*:H \ra W_y$ is an isomorphism, define $\partial H
: \wedge^2 W_y \ra T_xM/W_y$ by
$$\partial H (V_1, V_2) = [\tilde{V}_1, \tilde{V}_2] \mod W_y,$$
\noindent $\tilde{V}_i, i=1,2$ is a local vector field on $M$ with
$\tilde{V}_i(x') \in W_{\psi(x')}$ for a local section $\psi$ of
$\CF$ with $\psi_*(W_y)=H$. It is well defined(\cite{G}).

A $k$-subspace $H$ of $T_y\CF$
 is said to be a {\it 2-jet of an integral
variety} if $\pi_*$ restricts to an isomorphism $H \subset T_y\CF
\ra W_y \subset T_xM$ and $\partial H=0$.
\end{defn}

Such an $H$ is  indeed  a candidate for the 2-jets of actual
integral varieties of $\CF$ and the set of such an $H$ is again a
subvariety of $Gr(k, V+T_fF)$.

\begin{prop}
Let $\CF$ be an $F$-structure on $M$. Let $x \in M$ and let $y \in
\CF_x$. If $X$ is an integral variety of $\CF$ passing through $x
$ and tangent to $W_y$,  then the tangent space $H=T_y\tilde{X}
\subset T_y \CF$ of the lifting $\tilde{X} :=\{ (x, [T_xX]) | x
\in X, [T_xX] \in Gr(k, T_xM) \}$ of $X$ to $Gr(k, T_xM)$
satisfies

(1) $\pi_* : H \ra W_y$ is an isomorphism

(2) $\partial H=0$.

\end{prop}

\begin{pf} See Chapter 1 of \cite{G}. \end{pf}

The condition that $\pi_*: H \ra W_y$ is an isomorphism is
equivalent to the condition that $H$ is the graph of a map $p: W_y
\ra T_y(\CF_x) \subset W_y^* \otimes (T_xM/W_y)$. Define $\partial
p :\wedge^2 W_y \ra T_xM /W_y$ by $\partial p(V_1,
V_2)=p(V_1)(V_2)-p(V_2)(V_1)$, considering $T_y(\CF_x)$ as a
subspace of $W_y^* \otimes (T_xM/W_y)$. Then $\partial H=0$ if and
only if $\partial p=0$.

\begin{defn} \label{prolongation}
Let $\CF$ be an $F$-structure on $M$. For each $y \in \CF$, define
$\CF^{(1)}_y$ by the set of all 2-jets of integral varieties
tangent to $W_y$, i. e.,
$$\CF^{(1)}_y =\{ H \subset T_y\CF :  \pi_* : H \ra W_y \text{ is an
isomorphism and } \partial H=0 \}.$$ \noindent We call $
\CF^{(1)}=\cup_{y \in \CF} \CF^{(1)}_y $ the {\it first
prolongation} of $\CF$.
 Define
$$
F^{(1)}=\{ H \subset V+T_fF | H \text{ is the graph of a map } p :
W_f \ra T_fF,\,\,
\partial p =0 \}.$$
\noindent Then $\CF^{(1)}$ is an $F^{(1)}$-structure on $\CF$.
\end{defn}

Put $H^{1,1}(F)=Ker( \partial :W_f^* \otimes T_fF \ra \wedge^2
W_f^* \otimes (V/W_f))$. If $H^{1,1}(F)=0$ then $F^{(1)}$ is just
a point and the first prolongation defines a distribution on
$\CF$. In this case, the integrability of this distribution is
equivalent to the integrability of the $F$-structure $\CF$. So
there is at most one integral manifold passing through a given
point and tangent to a given $k$-subspace of the tangent space.
For the details see Chapter 1 of \cite{G}.



\section{Schubert rigidity}

\subsection{Description of the Schubert differential system}

\begin{exam}
For $\ba=(p^q)^*$, The Schubert variety $\sigma_{\ba}$ is the
sub-Grassmannian $\{ E \in Gr(m,n)| \BC^{m-q} \subset E \subset
\BC^{m+p} \} \simeq Gr(q, p+q)$. The Schubert differential system
$\CB_{\ba}$ is the flag space $F(m-q,m,m+p,n)$ and the parameter
space of the family of Schubert varieties of type $\ba$ is
$F(m-q,m+p,n)$. In this case, there is a double fibration $F(m-q,
m+p, n) \leftarrow \CB_{\ba} \rightarrow Gr(m,n)$.
\end{exam}

In general, the Schubert differential systems $\CB_{\ba}$ is a
generalized flag variety. Since $SL(n)$ acts on $\CB_{\ba}$
transitively and $\CB_{\ba}$ is compact, we have only to find the
corresponding subset of simple root system which generate the
isotropy group which is parabolic.

 Let $\CS=\{ \alpha_1, \cdots,
\alpha_{n-1}\}$ be the set of simple roots of $G=SL(n)$ and let
$P$ be the parabolic subgroup of $SL(n)$ generated by
$\CS^1=\CS-\{\alpha_m \}$, i.e. $\CS^1$ is the set of simple roots
of the semisimple part $SL(m) \times SL(n-m)$ of $P$. Then
$Gr(m,n)$ is $G/P$.

\begin{prop} \label{surj}
For $\ba=(p_1^{q_1}, \cdots, p_r^{q_r}), p_r \not=0$, put
$$\CS_{\ba}=\CS^1 -\{ \alpha_{q_1}, \cdots, \alpha_{q_1 +
\cdots+ q_r} \} \cup \{ \alpha_{n-p_1} \cdots, \alpha_{n-p_r} \}
.$$ 
Let $Q_{\ba}$ is the parabolic subgroup of $SL(n)$ generated by
$\CS_{\ba} \cup \{ \alpha_m \}$ and $P_{\ba}$ be the parabolic
subgroup  of $SL(m) \times SL(n-m)$ generated by the set
$\CS_{\ba}$. Then $\CB_{\ba}$ is the homogeneous manifold
$G/(Q_{\ba} \cap P)$ and $P_{\ba}$ is the isotropy group of the
the action of $SL(m) \times SL(n-m)$ on the fiber $B_{\ba}$ of
$\CB_{\ba} \ra Gr(m,n)$.
  If $\CS_{\ba} \subset \CS_{\bb}$, then there exists
a quotient map $\varphi_{\ba,\bb}: \CB_{\ba} \ra \CB_{\bb}$ which
preserves the fibers.
\end{prop}

For example, the following figure describes $\fn_{\ba}$,
$\fm_{\ba}$ and $Q_{\ba}$ for  $n=10, m=4$ and $\ba=(6,4,2,2)$,
where $\fm_{\ba}$ denotes the tangent space of $B_{\ba}$.
$\fn_{\ba}$ is the space of all $n \times n$-matrix with nonzero
elements only in $\ast$ and $\fm_{\ba}$ is the space of all $n
\times n$-matrix with nonzero elements only in $\bullet$. The Lie
algebra of the reductive part of $Q_{\ba}$ is the space of all $n
\times n$-matrix with nonzero elements only in $\diamond$.

\vs 

\begin{tabular}{|c|c|c|c||c|c|c|c|c|c|} \hline
 $\diamond$  & & & & & & & & & \\ \hline
$\bullet$ & $\diamond$ & &  & & & & & & \\ \hline

$\bullet$ & $\bullet$ & $\diamond$ & $\diamond$  & & & & && \\
\hline

$\bullet$ & $\bullet$ &$\diamond$ & $\diamond$ & & & & & & \\
\hline \hline

 & $\ast$ & $\ast$ & $\ast$ & $\diamond$ & $\diamond$ & & & & \\ \hline

 & $\ast$ & $\ast$ & $\ast$ & $\diamond$ & $\diamond$ & & & & \\ \hline

 &  & $\ast$ & $\ast$ & $\bullet$ & $\bullet$ & $\diamond$ & $\diamond$ & & \\ \hline

 &  & $\ast$ & $\ast$ & $\bullet$ & $\bullet$ & $\diamond$ & $\diamond$& & \\ \hline

  & & & &$\bullet$ & $\bullet$ & $\bullet$ & $\bullet$ & $\diamond$ & $\diamond$ \\
  \hline

    & & & &$\bullet$ & $\bullet$ & $\bullet$ & $\bullet$ & $\diamond$ & $\diamond$ \\ \hline
\end{tabular}

\vs

\begin{prop} \label{smooth}
For $\ba=(p^q)^*$,  the Schubert variety $\sigma_{\ba}$ of type
$\ba$ is strongly rigid except when ($p=1$ and $q \not=m$) or ($p
\not =n-m$ and $q=1$).
\end{prop}

\begin{pf} See \cite{W}. \end{pf}

\subsection{Foliation by smooth Schubert varieties}

We start with the simplest case and then will use the induction to
prove the general case.

\begin{exam} \label{exam}

The case when $\ba=(p^q)$ is studied in Example 13 and Remark 33
of \cite{B}.
 Fix a  $(n-m-p+q)$-subspace $\Lambda=\BC^{n-m-p+q}$ of
$\BC^n$. Then $\sigma_{\ba}(\Lambda)$ can be expressed as a union
of a family of  Schubert varieties of type
$\bb=((n-m)^q)=((n-m)^{(m-q)})^*$.

\begin{eqnarray*}
\sigma_{\ba}(\Lambda)&=& \{E \in Gr(m,n)| \dim(E \cap
\Lambda) \geq q \} \\
&=& \bigcup_{\BC^q \subset \Lambda} \{E \in Gr(m,n) | \BC^q
\subset E \}
\end{eqnarray*}


 \noindent Note that $Gr(q,n)$ is the parameter space of
 the Schubert varieties of type $\bb$ and $\{\BC^q \in Gr(q,n)| \BC^q \subset \Lambda \}$ is
 a Schubert variety of type
 $\bc=((n-m-p)^q)^*$
 in $Gr(q,n)$.

\vs

\vs $\ba=$\begin{tabular}{|c|c||c|c|}
\hline 0 &0& 0&0  \\
\hline $\fm_{\ba}$ &0&0&0 \\
\hline \hline
$\fn_{\ba}$ & $\fn_{\ba}$ & 0&  0 \\
\hline 0 & $\fn_{\ba}$ & $\fm_{\ba}$ & 0 \\
\hline
\end{tabular}
\quad $\bb=$\begin{tabular}{|c|c||c|c|}
\hline 0 &0& 0&0  \\
\hline $\fm_{\bb}$ &0&0&0 \\
\hline \hline
0 & $\fn_{\bb}$ & 0&  0 \\
\hline 0 & $\fn_{\bb}$ & 0 & 0 \\
\hline
\end{tabular}
\quad $\bc=$\begin{tabular}{|c||c|c|c|}
\hline 0 &0& 0&0  \\
\hline \hline 0 &$\fm_{\bc}$&0&0 \\
\hline
$\fn_{\bc}$ & 0 & 0&  0 \\
\hline 0 & $\fm_{\bc}$ & 0 & 0 \\
\hline
\end{tabular}
\vs

From this expression, we get the following desingularisation
$\pi_2$ of $\sigma_{\ba}(\Lambda)$.
$$\begin{array}{c}
 \{ (\BC^q, E) | \BC^q \subset E \}
 \stackrel{\pi_2}{\longrightarrow}  Gr(m,n) \supset
\sigma_{\ba}(\Lambda) \\
\pi_1 \, \downarrow \qquad \qquad \qquad \qquad  \qquad \qquad \qquad\\
 Gr(q,n) \supset Gr(q,\Lambda) \qquad \qquad  \qquad \qquad
\end{array}$$

\noindent Here, $\sigma_{\ba}(\Lambda)$ is equal to $
\pi_2(\pi_1^{-1}(Gr(q, \Lambda)))$ and $\pi_2 : \pi_1^{-1}(Gr(q,
\Lambda)) \ra \sigma_{\ba}(\Lambda)$ is generically one-to-one.

The smooth locus of $\sigma_{\ba}(\Lambda)$ is foliated by
Schubert varieties of type $\bb$ in $Gr(m,n)$ and the space of
leaves of this foliation is a Schubert variety of type $\bc$ in
$Gr(q,n)$.

\vs In the  same way as above, we get a different
desingularisation of $\sigma_{\ba}(\Lambda)$ by considering
$\sigma_{\ba}(\Lambda)$ as a union of a family of another type of
Schubert varieties:
$$\sigma_{\ba}(\Lambda)= \bigcup_{\Lambda \subset \BC^{n-p}} \{E
\in Gr(m,n) | E \subset \BC^{n-p} \} .$$

\end{exam}

\vs

\begin{prop} \label{foliation} For the partition $\ba=(p^q)$, the
Schubert variety $\sigma_{\ba}$ of type $\ba$ is rigid if $p >1$
and $q
>1$.
\end{prop}

It is proved in Example 13 and Remark 33 of  \cite{B}. We will
prove it again in such a way that  can be generalized to the case
of the Schubert differential systems of other singular Schubert
varieties.

\begin{lem} \label{foli} Let $\ba$ and $\bb$ be two partitions such
that $\fn_{\bb}$ is  a subspace of $\fn_{\ba}$ with $\CS_{\ba}
\subset \CS_{\bb}$ (and thus there is a projection $\varphi_{\ba,
\bb}:\CB_{\ba} \ra \CB_{\bb}$ as in proposition \ref{surj}).
Assume that $\CB_{\bb}$ is strongly rigid. For $p: \fn_{\ba} \ra
\fm_{\ba}$, define $\tilde{p}$ by the composition $ \fn_{\ba} \ra
\fm_{\ba} \twoheadrightarrow \fm_{\bb}$. If $\tilde{p}(v)=0$ for
all $v \in \fn_{\bb}$ and for all $p \in H^{1,1}(B_{\ba})$, then
the smooth locus of any integral variety of $\CB_{\ba}$ is
foliated by integral varieties of $\CB_{\bb}$, which are Schubert
varieties of type $\bb$.

\end{lem}

\begin{pf} By the strong rigidity of  $\CB_{\bb}$, its first
prolongation $\CB_{\bb}^{(1)}$ gives a distribution $D$ on
$\CB_{\bb}$ which is integrable because $\CB_{\bb}$ is integrable.
Integral varieties of $D$ are isomorphic to Schubert varieties of
type $\bb$ via the map $\pi_{\bb}:\CB_{\bb} \ra Gr(m,n)$.

Let $X$ be an integral variety of $\CB_{\ba}$. Let $\tilde{X}
\subset \CB_{\ba}$ be the lifting of $X$, that is,
$\tilde{X}=\{(x, [T_xX])| x \in X \}$. It suffices to show that
$\varphi_{\ba, \bb}(\tilde{X})$ is foliated by the integral
varieties of the distribution $D$ induced by $\CB^{(1)}_{\bb}$,
that is, at each point $y \in \tilde{X}$, $(\varphi_{\ba,
\bb})_*(T_y\tilde{X})$ contains $D_{\varphi_{\ba, \bb}(y)}$.

$$\begin{array}{c}
\qquad \quad \quad \CB_{\bb}^{(1)} \\
\,\,\,\qquad \downarrow \\
 \qquad   \tilde{X} \subset \CB_{\ba}
 \stackrel{\varphi_{\ba, \bb}}{\longrightarrow} \CB_{\bb} \supset
\varphi_{\ba, \bb}(\tilde{X})  \\
\searrow  \,\,\swarrow \,\, \\
X \subset Gr(m,n) \qquad
\end{array}
$$

The map $(\varphi_{\ba, \bb})_*$ is given by the projection
$\fm+\fm_{\ba} \ra \fm+\fm_{\bb}$ and $T_y\tilde{X}$ is the graph
of a map $p:\fn_{\ba} \ra \fm_{\ba}$ with $\partial p=0$. So
$(\varphi_{\ba, \bb})_*(T_y\tilde{X})$ is the graph of the map
$\fn_{\ba} \stackrel{p}{\ra} \fm_{\ba} \ra \fm_{\bb}$ with
$\partial p=0$. On the other hand $D_{\varphi_{\ba,\bb}(y)}$ is
the graph of the zero map $\fn_{\bb} \ra \fm_{\bb}$. By the
assumption, $\tilde{p}(v)$ is zero for all $v \in \fn_{\bb}$, so
$\varphi_{\ba, \bb}(T_y\tilde{X})$ contains $D_{\varphi_{\ba,
\bb}(y)}$.
\end{pf}

\begin{lem}  \label{A2} Let $\ba, \bc$ be as in Example \ref{exam}.
Let $A \subset Gr(q,n)$ be a subvariety of dimension $q(n-m-p)$.
Define $X_A:= \pi_2(\pi_1^{-1}(A)) \subset Gr(m,n)$.  If $X_A$ is
an integral variety of $\CB_{\ba}(Gr(m,n))$ and $ dim
(\pi_1^{-1}(A))$ is equal to $dim(X_A)$, then $A$ is an integral
variety of $\CB_{\bc}(Gr(q,n))$.

$$\begin{array}{c}
 \qquad  \qquad   \CB_{\ba} \ra \CB_{\bb}
 \stackrel{\pi_1}{\longrightarrow} Gr(q,n) \supset A  \\
\searrow  \,\,\swarrow  \pi_2 \qquad  \\
\,\, Gr(m,n) \,\,\, \supset X_A
\end{array}
$$

\end{lem}

\begin{pf} We will follow the arguments in Example 16.6 of \cite{H}.
 Define $\tilde{\pi_i}, i=1,2$ to be the projection
 from $Gr(q, n) \times Gr(m,n)$ to the first
and second component, respectively. Let $F(q,m,n)$ denote the flag
space $\{(\Gamma, E) | \Gamma \subset E \subset \BC^n \} \subset
Gr(q, n) \times Gr(m,n) $. Then $\pi_2(\pi_1^{-1}(A))$ is
$\tilde{\pi}_2(\tilde{\pi}_1^{-1}(A) \cap F(q,m,n))$.

Let $\Gamma$ be a smooth point in $A$. Then
$\tilde{\pi}_1^{-1}(A)$ is smooth at $(\Gamma, E)$ for all $E \in
Gr(m,n)$ with the tangent space
$$T_{(\Gamma, E)}\tilde{\pi}_1^{-1}(A) = \left\{(\eta, \varphi)|
\begin{array}{ll} \eta : \Gamma \ra \BC^n/\Gamma, & \eta \in
T_{\Gamma}A \\
                  \varphi : E \ra \BC^n/E, &
                  \end{array}
                   \right\}. $$
 The tangent space
of $F(q,m,n)$ at $(\Gamma, E)$ is
$$T_{(\Gamma, E)}F(q,m,n) = \left\{(\eta, \varphi)|
\begin{array}{ll} \eta : \Gamma \ra \BC^n/\Gamma, &  \\
                  \varphi : E \ra \BC^n/E, & \varphi|_{\Gamma} \equiv
\eta \mod E
                  \end{array}
                   \right\}. $$
By dimension counting, we see that the two tangent spaces are
transversal so that $\pi_1^{-1}(A)=\tilde{\pi}_1^{-1}(A) \cap
F(q,m,n)$ is smooth at all $(\Gamma, E) $ with $\Gamma \subset E$,
and the tangent space of $\pi_1^{-1}(A)$ at $(\Gamma, E)$
 is given by
$$T_{(\Gamma, E)}\pi_1^{-1}(A) = \left\{(\eta, \varphi)|
\begin{array}{ll} \eta : \Gamma \ra \BC^n/\Gamma, &  \eta \in
T_{\Gamma}A\\
                  \varphi : E \ra \BC^n/E, & \varphi|_{\Gamma} \equiv
\eta \mod E
                  \end{array}
                   \right\}. $$

For each $E \in X_A$ there are only finitely many $\Gamma \in A$
with $\Gamma \subset E$. If there are more than one $\Gamma \in A$
with $\Gamma \subset E$, then $X_A$ is not smooth at
$E$(Proposition 16.8 of \cite{H}). Let $E$ be an element in $X_A$
such that there is only one $\Gamma \in A$ with $\Gamma \subset
E$. Then $\pi_2: \pi_1^{-1}(A) \ra X_A$ is one-to-one over $E$ so
$X_A$ is smooth at $E$ with the tangent space
 $$T_E(X_A)=\{ \varphi : E \ra \BC^n/E \,\,| \,\,\varphi|_{\Gamma}
\equiv \eta \mod E,\,\, \eta \in T_{\Gamma}A,\}.$$ \noindent Since
$T_E(X_A)$ is of type $\ba$ in $Gr(m,n)$, $T_{\Gamma}A$ is of type
$\bc$ in $Gr(q,n)$.
\end{pf}

{\bf Proof of Proposition \ref{foliation}}
We will show that for all $p:\fn_{\ba} \ra \fm_{\ba}$ with $p \in
H^{1,1}(B_{\ba})$, the induced map $\tilde{p}: \fn_{\bb} \subset
\fn_{\ba} \ra \fm_{\ba} \ra \fm_{\bb}$  is zero. Then for an
integral variety $X$ of $\CB_{\ba}$, the space of leaves
$A=\pi_1(\varphi_{\ba, \bb}(\tilde{X})) \subset Gr(q,n)$ of the
foliation on $X$ given by Lemma \ref{foli} satisfies the
conditions in Lemma \ref{A2}. Thus $A \subset Gr(q,n)$ is an
integral variety of $\CB_{\bc}(Gr(q,n))$. Since $q
>1$ $\CB_{\bc}$ is rigid. So $A$ is the Schubert variety
$Gr(q,\Lambda)$ for a $(n-m-p+q)$-subspace $\Lambda$ of $\BC^n$
and hence $X$ is the Schubert variety $\sigma_{\ba}(\Lambda)$.

\vs Let $E_{i,j}$ be the $n \times n$-matrix with only one nonzero
element in the $i$-th row and $j$-th column. Then $[E_{i,j},
E_{k,\ell}]=\delta_{j,k} E_{i \ell} -\delta_{i, \ell} E_{k,j}$ for
all $1 \leq i,j, k, \ell \leq n$.

\vs Since $[p(X), Y]-[p(Y), X] \in \fn_{\ba}$ for all $X, Y \in
\fn_{\ba}$ and $[\fm_{\ba}/\fm_{\bb}, \fn_{\bb}] \subset
\fn_{\bb}$, $[\tilde{p}(X), Y]-[\tilde{p}(Y), X] \in \fn_{\ba}$
for all $X, Y \in \fn_{\bb}$.  Assume that $q+1 \leq i,j \leq m$
and $m+1 \leq r,s \leq n-p $ and $n-p+1 \leq a,b \leq n$. Note
that for a fixed $a$,   if $X \in \fm_{\bb}$ and  $[X, E_{a,i}]=0$
for all $i$, then $X=0$.


\vs
\begin{tabular}{|c|c||c|c|}
\hline 0 &0& 0&0  \\
\hline $\fm_{\bb}$ &0&0&0 \\
\hline \hline
$\ast$ & $E_{r,i}$ & 0&  0 \\
\hline 0 & $ E_{a, i}$ & \,$\bullet$\, & \,0\, \\
\hline
\end{tabular}

\vs


 From $[\tilde{p}(E_{a,i}),
E_{r,j}]-[\tilde{p}(E_{r,j}), E_{a,i}] \in \fn_{\ba}$, we get both
$[\tilde{p}(E_{a,i}), E_{r,j}]$ and $[\tilde{p}(E_{r,j}), E_{a,i}]
$ are contained in $\fn_{\ba}$. Thus $[\tilde{p}(E_{r,j}),
E_{a,i}] $ is zero, which implies that $\tilde{p}(E_{r,j})=0$.

Since $[\tilde{p}(E_{a,i}), E_{b,j}]-[\tilde{p}(E_{b,j}),
E_{a,i}]$ is contained in $\fn_{\ba}$, it should be zero and thus
both $[\tilde{p}(E_{a,i}), E_{b,j}]$ and $[\tilde{p}(E_{b,j}),
E_{a,i}]$ should be zero for $a \not=b$. Here we use the condition
that $p
> 1$. So $\tilde{p}(E_{b,j})$ is zero.
 $\square$

\subsection{Sub-Grassmannians}

 To extend the Proposition \ref{foliation} to the general case,
we consider the following problem: Suppose that a Schubert variety
$\sigma_{\ba}$ is contained in a proper sub-Grassmannian of
$Gr(m,n)$ and $\sigma_{\bb}$ is the minimal sub-Grassmannian among
them. Then will any integral variety of $\CB_{\ba}$ be contained
in a sub-Grassmannian $\sigma_{\bb}$?

\begin{prop}  \label{subgr} Let $\ba=(p_1^{q_1}, \cdots, p_r^{q_r}) \in
P(m,n)$ be a partition and let $\ba'=({p'_1}^{q'_1}, \cdots,
{p'_r}^{q'_r})$ be the conjugate of $\ba$. Suppose that
$\fn_{\ba}$ is contained in a proper rectangle in $\fm$. Let $\bb$
be the partition corresponding to the minimal rectangle among
them. Then any integral variety of $\CB_{\ba}$ is contained in
  a sub-Grassmannian $\sigma_{\bb}$ except when
 $q_1 + \cdots + q_r=m$ and  $q_r =1$ or
$q'_1 + \cdots + q'_r=n-m$ and $q'_r =1$.

\end{prop}

If both $q_1 + \cdots + q_r < m$ and $q'_1 + \cdots + q'_r<n-m$
hold, then there is no proper sub-Grassmnnian containing
$\sigma_{\ba}$. So the cases we will consider below is either when
$q_1 + \cdots + q_r=m$ and  $q_r \geq 2$ or when $q'_1 + \cdots +
q'_r=n-m$ and $q'_r \geq 2$.

\begin{lem} \label{inclusion} Let $\ba$ and $\bb$ be two
partitions such that $\fn_{\ba}$ is a subspace of $\fn_{\bb}$ with
$\CS_{\ba} \subset \CS_{\bb}$ and thus there is a projection
$\CB_{\ba} \ra \CB_{\bb}$ as in proposition \ref{surj}. Assume
that $\CB_{\bb}$ is strongly rigid. For $p: \fn_{\ba} \ra
\fm_{\ba}$, define $\tilde{p}$ by the composite map $ \fn_{\ba}
\ra \fm_{\ba} \twoheadrightarrow \fm_{\bb}$. If $\tilde{p}=0$ for
all $p \in H^{1,1}(B_{\ba})$, then any integral subvariety of
$\CB_{\ba}$ is contained in an integral variety of $\CB_{\bb}$,
which is a Schubert variety of type $\bb$.
\end{lem}

\begin{pf} The proof is similar to the proof of Lemma \ref{foli}.
\end{pf}

\vs \noindent {\bf Proof of Proposition \ref{subgr}}
%
%
%
%
First, we  consider the case  $\ba=(p_1^{q_1}, p_2^{q_2},
p_3^{q_3})$, $p_3 \not=0$ is a partition with $q_1 + q_2 +q_3=m$
and $q_3 \geq 2$. As the proof will show, the general case can be
obtained in the same way.

\vs {\bf Claim.}  Put $\bb=(p_3^{m})$.

\vs $\ba=$
\begin{tabular}{|r|r|r||r|r|r|}
\hline 0 & 0&0&0&0&0 \\
\hline $\bullet$&0&0&0&0&0 \\
\hline $\bullet$&$\bullet$ & 0&0&0&0 \\
\hline \hline 0 & $\ast$ & $\ast$ & 0&0&0 \\
\hline 0 & 0 & $\ast$ &$ \bullet$ &0&0 \\
\hline 0 & 0 & 0 &$\bullet$&$\bullet$&0\\
\hline
\end{tabular}
\quad $\bb=$
\begin{tabular}{|r|r|r||r|r|r|}
\hline 0&0&0&0&0&0 \\
\hline 0&0&0&0&0&0 \\
\hline 0&0&0&0&0&0 \\
\hline \hline $\ast$ & $\ast$ & $\ast$ &0&0&0\\
\hline $\ast$ & $\ast$ & $\ast$ &0&0&0\\
\hline 0 & 0 & 0 &$\bullet$&$\bullet$&0\\
\hline
\end{tabular}

\vs

\noindent Then any integral variety of $\CB_{\ba}$ is contained in
a Schubert variety of type $\bb$.

\vs

 {\bf Proof.} Let $p: \fn_{\ba} \ra \fm_{\ba}$ be a map with $\partial p=0$
 Then $[p(X), Y]-[p(Y),X] \in \fn_{\ba}$ for all $X, Y \in \fn_{\ba}$.
Put $\tilde{p} : \fn_{\ba} \ra \fm_{\bb}$ to be the composition of
$p$ with the projection $\fm_{\ba} \ra \fm_{\bb}$.

Since  $[\fm_{\ba} / \fm_{\bb}, \fn_{\ba}] \subset \fn_{\bb}$, we
have $[\tilde{p}(X), Y] -[\tilde{p}(Y), X] \in \fn_{\bb}$ for all
$X, Y \in \fn_{\ba}.$ Assume that $q_1+1 \leq k, \ell \leq q_1 +
q_2$, $q_1+q_2 +1 \leq i,j \leq m$, $m+1 \leq r,s \leq n-p_1$ and
$n-p_1 +1 \leq a,b \leq n-p_2$. Note that for $X \in \fm_{\bb}$,
if $[X, E_{r,i}]=[X, E_{a,i}]=0$ for all $r$ and $a$ for a fixed
$i$, then $X$=0.

\vs
\begin{tabular}{|c|c|c||c|c|c|}
\hline \,\,0 \,& 0&0&0&0&\,\,0\,\, \\
\hline $\bullet$&0&0&0&0&0 \\
\hline $\bullet$&$\bullet$ & 0&0&0&0 \\
\hline \hline 0 & $E_{r,k}$ & $E_{r,i}$ & 0&0&0 \\
\hline 0 & 0 & $E_{a,i} $ &$ \bullet$ &0&0 \\
\hline 0 & 0 & 0 &$\fm_{\bb}$&$\fm_{\bb}$&0\\
\hline
\end{tabular}
\vs

\noindent From $[\tilde{p}(E_{r,k}), E_{s,i}]-[\tilde{p}(E_{s,i}),
E_{r,k}] \in \fn_{\bb}$, we see that both
 $[\tilde{p}(E_{r,k}), E_{s,i}]$ and
 $[\tilde{p}(E_{s,i}), E_{r,k}] $ are contained in
 $\fn_{\bb}$. So $[\tilde{p}(E_{r,k}), E_{s,i}]=0$.
 The same equation hold if we replace $E_{s,i}$ by $E_{a,i}$.
 Thus $\tilde{p}(E_{r,k})=0$.

Put $\bc=((n-m)^{(q_1+q_2)}, p_3^{q_3})$.

\vs $\bc=$
\begin{tabular}{|r|r|r||r|r|r|}
\hline 0&0&0&0&0&0 \\
\hline 0&0&0&0&0&0 \\
\hline $\bullet$&$\bullet$&0&0&0&0 \\
\hline \hline 0 & 0 & $\ast$ &0&0&0\\
\hline 0 & 0 & $\ast$ &0&0&0\\
\hline 0 & 0 & 0 &$\bullet$&$\bullet$&0\\
\hline
\end{tabular}
\vs

\noindent Then $\tilde{p}$ restricts to a map $\fn_{\bc} \ra
\fm_{\bb} \subset \fm_{\bc}$. Also we have $[\tilde{p}(X), Y]
-[\tilde{p}(Y), X] \in \fn_{\bc}$ for all $X, Y \in \fn_{\bc}$.
Since $q_3 \geq 2$, $\CB_{\bc}$ is strongly rigid so that
$\tilde{p}$ is zero on $\fn_{\bc}$. Hence $\tilde{p}$ is zero. By
Lemma \ref{inclusion}, any integral variety of $\CB_{\ba}$ is
contained in a Schubert variety of type $\bb$.

\vs In the same way, we can show that any integral manifold of
$\CB_{\ba}$ is contained in a Schubert variety of type $\hat{\bb}$
for $\hat{\bb}=({p'_3}^{n-m})$.


\vs$\ba=$
\begin{tabular}{|r|r|r||r|r|r|}
\hline 0 & 0&0&0&0&0 \\
\hline $\bullet$&0&0&0&0&0 \\
\hline $\bullet$&$\bullet$ & 0&0&0&0 \\
\hline \hline 0 & $\ast$ & $\ast$ & 0&0&0 \\
\hline 0 & 0 & $\ast$ &$ \bullet$ &0&0 \\
\hline 0 & 0 & 0 &$\bullet$&$\bullet$&0\\
\hline
\end{tabular}
\quad $\hat{\bb}=$
\begin{tabular}{|r|r|r||r|r|r|}
\hline 0&0&0&0&0&0 \\
\hline $\bullet$&0&0&0&0&0 \\
\hline $\bullet$&0&0&0&0&0 \\
\hline \hline 0 & $\ast$ & $\ast$ &0&0&0\\
\hline 0& $\ast$ & $\ast$ &0&0&0\\
\hline 0 & $\ast$ & $\ast$ &0&0&0\\
\hline
\end{tabular}

\vs

\noindent Here we use the rigidity of the Schubert differential
system $\CB_{\hat{\bc}}$ with $\hat{\bc}=( m^{(q'_1+q'_2)},
{p'_3}^{q'_3}), q_3' \geq 2$.

\vs  $\hat{\bc}=$
\begin{tabular}{|r|r|r||r|r|r|}
\hline 0&0&0&0&0&0 \\
\hline $\bullet$&0&0&0&0&0 \\
\hline $\bullet$&0&0&0&0&0 \\
\hline \hline 0 & $\ast$ & $\ast$ &0&0&0\\
\hline 0 & 0 &0 &$\bullet$&0&0\\
\hline 0 & 0 & 0 &$\bullet$&0&0\\
\hline
\end{tabular}

$\square$

\begin{thm} \label{Schubertgr}
Let $\ba=(p_1^{q_1}, \cdots, p_r^{q_r})$ be a partition and let
$\ba'=({p'_1}^{q'_1}, \cdots, {p'_r}^{q'_r})$ be its conjugate.
Then $\sigma_{\ba}$ is Schubert rigid if
 $q_i, q_i' \geq 2$ for all $i \leq r$.

\end{thm}

\begin{pf} We will use the induction on $r$. Thanks to Proposition
\ref{subgr}, we may assume that there is no proper sub-Gassmannian
containing $\sigma_{\ba}$ so that $q_1 + \cdots q_r <m$ and $q'_1
+ \cdots q'_r <n-m$.

By Lemma \ref{foli} and the same argument as in Proposition
\ref{foliation}, any integral varieties of $\CB_{\ba}$ is foliated
by the Schubert varieties of type $\bb=((n-m)^q)$, where $q=q_1+
\cdots q_r$. Then the space of leaves $A$ of this foliation will
be an integral variety of the Schubert differential system
$\CB_{\bc}$ on $Gr(q, n)$, $\bc=((p_1+(m-q))^{q_1}, \cdots, (p_r +
(m-q))^{q_r} )$. By Proposition \ref{subgr}, any integral
varieties of $\CB_{\bc}$ is contained in a sub-Grassmannian of
type $((p_r+(m-q))^{q})$ in $Gr(q,n)$. Thus $A$ is an integral
varieties of the Schubert differential system $\CB_{\bf d}$, ${\bf
d}=((p_1-p_r)^{q_1}, \cdots, (p_{r-1}-p_r)^{q_{r-1}} )$ on this
sub-Grassmannian, which is rigid by the induction hypothesis.
\end{pf}

\section{Schur rigidity}

\subsection{Criterions for the equality $B_{\ba}=R_{\ba^*}$}

Let $\sigma_{\ba}$ be a Schubert variety. If $\sigma_{\ba}$ is
Schubert rigid and the Schur differential system $\CR_{\ba^*}$ is
equal to the Schubert differential system $\CB_{\ba}$, then
$\sigma_{\ba}$ is Schur rigid. To compare $\CB_{\ba}$ and
$\CR_{\ba^*}$, we will use the following description of
$\CR_{\ba^*}$.

\begin{prop} \label{hwo}
Let $B_{\ba}$ (resp. $R_{\ba^*}$) be the fiber of $\CB_{\ba}$
(resp. $\CR_{\ba^*}$) at $E \in Gr(m,n)$. Then $R_{\ba^*}$  is
equal to
$$Gr(|\ba^*|, E^* \otimes Q) \cap \PP(\BS_{\ba^*}(E^*) \otimes
\BS_{(\ba^*)'}(Q)) \subset \PP(\wedge^{|\ba^*|}(E^* \otimes Q))$$
and
 $B_{\ba}$ is the  orbit of a highest weight vector in the
irreducible representation space $\PP(\BS_{\ba^*}(E^*) \otimes
\BS_{(\ba^*)'}(Q))$ of $SL(E^*) \times SL(Q)$.

\end{prop}

\begin{pf}
By Proposition 3.1 of \cite{Ho}.
\end{pf}

Note that we adapt the convention that $\sigma_{\ba}$ is a
Schubert variety of dimension $|\ba^*|$, while  $X_{w}$ is a
Schubert variety of dimension $\ell(w)$ in \cite{Ho}. Thus, when
$\ba^* \in P(m,n)$ corresponds to $w \in W^P$, $\CB_{\ba}$ is
equal to $\CB_{w}$ and $\CR_{\ba^*}$ is equal to $\CR_w$.

\vs Since $P_0$ acts on the fiber $B_{\ba}$, we have a
decomposition $\fp_0=\fm_{\ba} + \fl_{\ba} + \fm_{\ba}^*$, where
the tangent space of $B_{\ba}$ at a point is isomorphic to
$\fm_{\ba}$. A sufficient condition for the equality
$B_{\ba}=R_{\ba^*}$ is given by Proposition 3.2 in \cite{Ho}. Put
$\bI_{\ba}=\BS_{\ba^*}(E^*) \otimes \BS_{(\ba^*)'}( Q)$.

\begin{prop} \label{sufficient}
 Assume that for the
$\fl_{\ba}$-generators $\varphi$ of all irreducible
$\fl_{\ba}$-representation spaces in the complement of $\fm_{\ba}$
in $\fn_{\ba}^* \otimes \fm/\fn_{\ba}$, we have
$$\varphi^k(v_1 \wedge \cdots \wedge v_k) \not \in \bI_{\ba},$$ \noindent
where  $\{v_1, \cdots, v_k \}$ is a basis of $\fn_{\ba}$ and
$\varphi^k: \wedge^k \fn_{\ba} \ra \wedge^k \fm/\wedge^k\fn_{\ba}$
is defined by $$\varphi^k(v_1 \wedge \cdots \wedge v_k)=\sum_i v_1
\wedge \cdots \wedge \varphi(v_i) \wedge \cdots \wedge v_k \mod
\wedge^k \fn_{\ba} .$$
 Then $B_{\ba}$ is equal to
$R_{\ba^*}$.
\end{prop}

Let $E$ be an $m$-subspace of $\BC^n$ and $Q$ be the quotient
$\BC^n /E$. Then $\fm$ is equal to $E^* \otimes Q$. Take a
partition $\ba =(p_1^{q_1}, \cdots, p_r^{q_r}) \in P(m,n)$ and let
$\ba'=({p'_1}^{q'_1}, \cdots, {p'_r}^{q'_r})$, $ p'_r \not=0$ be
its conjugate.  Write $E=\oplus_{i=1}^{r_E} E_i$ and
$Q=\oplus_{a=1}^{r_Q} Q_a$ so that $\fl_{\ba}=(\oplus_{i=1}^{r_E}
sl(E_i)) \oplus ( \oplus_{a=1}^{r_Q} sl(Q_a))$. Note that $r_E$ is
$r$ if $q_1 + \cdots +q_r=m$ and is $r+1$, otherwise, and $r_Q$ is
$r$ if $q_1' + \cdots q_r'=n-m$ and is $r+1$, otherwise. We
indexed $E_i$ and $Q_{a}$ keeping the order of the basis $\{ e_1,
\cdots, e_n\}$ of $\BC^n$ such that $\{ e_1, \cdots, e_m\}$ is a
basis of $E$. Set $r_i=\dim E_i$ and $s_a=\dim Q_a$ for $1 \leq i
\leq r_E$ and $1 \leq a \leq r_Q$. Then $E_i=<e_{r_{i-1}+1},
\cdots, e_{r_{i}}>$ and $Q_a=<q_{s_{a-1}+1}, \cdots, q_{s_a}>$,
where $q_p:=e_{m+p}$ for $p=1, \cdots, n-m$.


Let $\Pi$ be the index set of $(i,a)$ such that
$\fn_{\ba}=\oplus_{(i,a) \in \Pi } (E_i^* \otimes Q_a)$. As a
subspace of $\fn_{\ba}^* \otimes \fm/\fn_{\ba}$,
$\fm_{\ba}=(\oplus_{i <j} E_i^*\otimes E_j) \oplus (\oplus_{b <a}
Q_b^* \otimes Q_a)$ is equal to
$$(\oplus_{i <j} E_i^* \otimes E_j \otimes < Id_{\oplus_{a \in
\Pi_{i,j}}Q_a}>_{\BC}) \oplus ( \oplus_{b<a } < Id_{\oplus_{ i \in
\Pi_{b,a}}E_i}>_{\BC} \otimes Q_b^* \otimes Q_a)$$ \noindent where
$\Pi_{i,j}=\{a : (i,a) \not\in \Pi, (j,a) \in \Pi \}$ and
$\Pi_{b,a}=\{ i : (i,b) \in \Pi, (i,a) \not\in \Pi \}$.

We may choose the order of the set of roots of $SL(n)$ in such a
way that the maximal root is located in the most left and the
lowest box $E_1^* \otimes Q_s$ and the minimal root is located in
the most right and the highest box $E_r^* \otimes Q_1$. Then the
highest weight vector in $\fm_{\ba} \subset \fn_{\ba}^* \otimes
\fm/\fn_{\ba}$ is either

 $ \sum_{a \in
\Pi_{i,j}} \sum_{q_p \in  Q_a} x_{\alpha_p}^* \otimes
x_{\beta_p}$,
 where $x_{\alpha_p}=e_{r_j }^* \otimes q_p$ and $x_{\beta_p}=e_{r_{i-1}+1 }^*
 \otimes q_p$ for some $i <j$,
or,

 $ \sum_{i \in
\Pi_{b,a}} \sum_{e_p \in  E_i} x_{\alpha_p}^* \otimes
x_{\beta_p}$,
 where $x_{\alpha_p}=e_{p}^* \otimes q_{s_{b-1}+1}$ and $x_{\beta_p}=e_{p }^*
 \otimes q_{s_a}$ for some $b <a$,

\vs For example, Consider $\ba=(9^2, 7^2,3^4)$ in $P(10, 19)$.
Then we have $E=\oplus_{i=1}^4 E_i$, $Q=\oplus_{a=1}^3Q_a$ and
$\Pi=\{(2,1), (3,1),  (3,2), (4,1),  (4,2), (4,3) \}$. The highest
weight vector of $\fm_{\ba}$ is

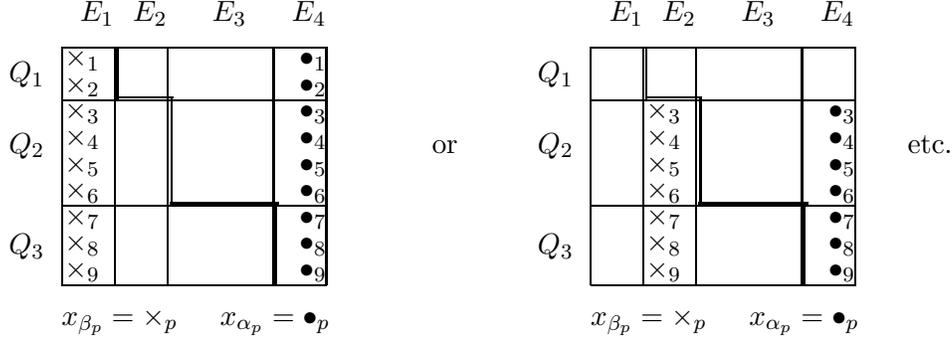
\begin{picture}(300, 130)

\put(0,0){\line(1,0){100}} \put(0,30){\line(1,0){100}}
\put(0,70){\line(1,0){100}} \put(0,90){\line(1,0){100}}

\put(0,0){\line(0,1){90}} \put(20,0){\line(0,1){90}}
\put(40,0){\line(0, 1){90}} \put(80,0){\line(0,1){90}}
\put(100,0){\line(0,1){90}}

\put(7, 100){$E_1$} \put(27, 100){$E_2$} \put(57, 100){$E_3$}
\put(87, 100){$E_4$}

\put (-20, 77){$Q_1$} \put(-20, 50){$Q_2$} \put(-20, 12){$Q_3$}


 \put(90, 83){$\bullet_1$}
 \put(90, 73){$\bullet_2$} \put(90, 63){$\bullet_3$}
 \put(90, 53){$\bullet_4$} \put(90, 43){$\bullet_5$}
 \put(90, 33){$\bullet_6$} \put(90, 23){$\bullet_7$}
 \put(90, 13){$\bullet_8$} \put(90, 3){$\bullet_{9}$}

 \put(1, 83){$\times_1$}
 \put(1, 73){$\times_2$} \put(1, 63){$\times_3$}
 \put(1, 53){$\times_4$} \put(1, 43){$\times_5$}
 \put(1, 33){$\times_6$} \put(1, 23){$\times_7$}
 \put(1, 13){$\times_8$} \put(1, 3){$\times_{9}$}



\put(200,0){\line(1,0){100}} \put(200,30){\line(1,0){100}}
\put(200,70){\line(1,0){100}} \put(200,90){\line(1,0){100}}

\put(200,0){\line(0,1){90}} \put(220,0){\line(0,1){90}}
\put(240,0){\line(0, 1){90}} \put(280,0){\line(0,1){90}}
\put(300,0){\line(0,1){90}}

\put(207, 100){$E_1$} \put(227, 100){$E_2$} \put(257, 100){$E_3$}
\put(287, 100){$E_4$}

\put(180, 77){$Q_1$} \put(180, 50){$Q_2$} \put(180, 12){$Q_3$}


 \put(290, 63){$\bullet_3$}
 \put(290, 53){$\bullet_4$} \put(290, 43){$\bullet_5$}
 \put(290, 33){$\bullet_6$} \put(290, 23){$\bullet_7$}
 \put(290, 13){$\bullet_8$} \put(290, 3){$\bullet_9$}

  \put(221, 63){$\times_3$}
 \put(221, 53){$\times_4$} \put(221, 43){$\times_5$}
 \put(221, 33){$\times_6$} \put(221, 23){$\times_7$}
 \put(221, 13){$\times_8$} \put(221, 3){$\times_9$}


\put(21, 90){\line(0,-1){20}} \put(21, 71){\line(1,0){21}}

\put(41.5, 71){\line(0,-1){41}} \put(41, 31){\line(1,0){41}}

\put(81, 31){\line(0,-1){31}}

\put(221, 90){\line(0,-1){20}} \put(221, 71){\line(1,0){21}}

\put(241.5, 71){\line(0,-1){41}} \put(241, 31){\line(1,0){41}}

\put(281, 31){\line(0,-1){31}}


\put(60, -15){$x_{\alpha_p}=\bullet_p$} \put(0,
-15){$x_{\beta_p}=\times_p$}

\put(260, -15){$x_{\alpha_p}=\bullet_p$} \put(200,
-15){$x_{\beta_p}=\times_p$}

\put(140, 50){or} \put(320, 50){etc.}
\end{picture}

\vs \vs \vs

\begin{prop} \label{typehwv} Let $E=\oplus_i E_i$ and $Q=\oplus_a Q_a$
be the decomposition associate to $\ba$ as in the above and $\{
e_1, \cdots, e_n \}$ be a basis for $\BC^n$ indexed as in the
above. Then the highest weight vector of an irreducible
$\fl_{\ba}$-representation space in the complement of $\fm_{\ba}$
in $\fn_{\ba}^* \otimes \fm /\fn_{\ba}$ is either

(1) a decomposable vector $x_{\alpha}^* \otimes x_{\beta}$, where
$x_{\alpha}=e_{r_j}^* \otimes q_{s_{b-1}+1}$ is the lowest weight
vector of $E_j^* \otimes Q_b$ and $x_{\beta}=e_{r_{i-1}+1}^*
\otimes q_{s_{a}}$ is the highest weight vector of $E_i^* \otimes
Q_a$ for some $(j,b) \in \Pi, (i,a) \not\in \Pi$, or,

 (2) $\sum_{q_p \in Q_a} x_{\alpha_p}^* \otimes x_{\beta_p}$,
 where $x_{\alpha_p}=e_{r_j }^* \otimes q_p$ and $x_{\beta_p}=e_{r_{i-1}+1 }^*
 \otimes q_p$ for some $(j,a) \in \Pi, (i,a) \not\in \Pi$
 such that $(j, a-1) \in \Pi$ and $(i, a-1) \not \in \Pi$, or,

 (3) $\sum_{e_p \in E_i} x_{\alpha_p}^* \otimes x_{\beta_p}$,
 where $x_{\alpha_p}=e_{p }^* \otimes q_{s_{b-1}+1}$ and $x_{\beta_p}=e_{p }^*
 \otimes q_{s_a}$ for some $(i,b) \in \Pi, (i,a) \not\in \Pi$
 such that $(i+1, b) \in \Pi$ and $(i+1, a) \not \in \Pi$.

\end{prop}

\begin{pf} If $i \not=j$ and $a \not=b$, then $(E_j^* \otimes Q_b)^*
\otimes (E_i^* \otimes Q_a)$ is an irreducible
$\fl_{\ba}$-representation space. But if $i \not=j$ and $a =b$,
then $(E_j^* \otimes Q_a)^* \otimes (E_i^* \otimes Q_a) \simeq E_j
\otimes E_i^* \otimes (Q_a^* \otimes Q_a)$ is decomposed as $(E_j
\otimes E_i^* \otimes <Id_{Q_a}>_{\BC} ) \oplus (E_j \otimes E_i^*
\otimes sl(Q_a))$, each of which are irreducible $sl(E_j) \times
sl(E_i) \times sl(Q_a)$-representation spaces. The highest weight
vector  of the irreducible representation spaces $(E_j^* \otimes
Q_b)^* \otimes (E_i^* \otimes Q_a)$ or $(E_j \otimes E_i^* \otimes
sl(Q_a))$ is of type (1).

The component $E_i^* \otimes E_j$ of $\fm_{\ba}$ corresponds to
the component $E_i^* \otimes E_j \otimes < Id_{\oplus_{a \in
\Pi_{i,j}}Q_a}>_{\BC}$ in $\fn_{\ba}^* \otimes \fm/\fn_{\ba}$, so
its complement  in $\oplus_{a \in \Pi_{i,j}} E_i^* \otimes E_j
\otimes <Id_{Q_a}>_{\BC}$ is $\oplus_{a \in \tilde{\Pi}_{i,j} }
E_i^* \otimes E_j \otimes < Id_{Q_a}>_{\BC}$, where
$\tilde{\Pi}_{i,j}=\{a : (i,a) \not\in \Pi, (j,a) \in \Pi, (i,a-1)
\not\in \Pi, (j, a-1) \in \Pi \}$  is obtained from $\Pi_{i,j}$ by
excluding the smallest index in $\Pi_{i,j}$. The highest weight of
the irreducible representation space in these components is of
type (2).

Considering the case when $i =j$ and $a \not=b$, we obtain the
highest weight vectors of type (3).
\end{pf}

For example, the highest weight vector $\sum_{q_p \in Q_a}
x_{\alpha_p}^* \otimes x_{\beta_p}$ of type (2) is

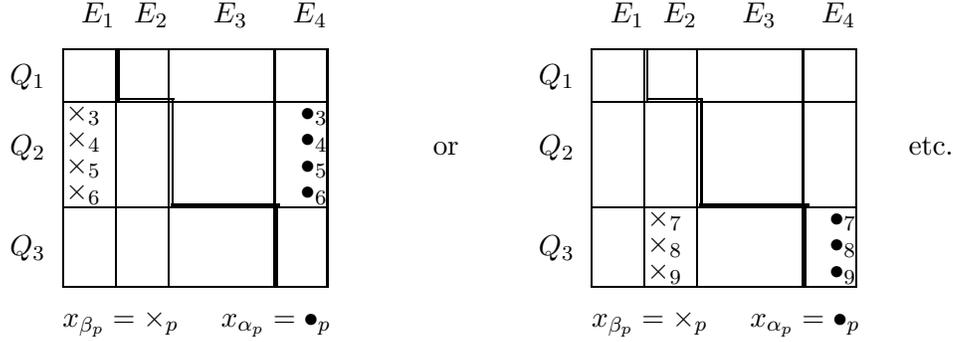
\begin{picture}(300, 130)

\put(0,0){\line(1,0){100}} \put(0,30){\line(1,0){100}}
\put(0,70){\line(1,0){100}} \put(0,90){\line(1,0){100}}

\put(0,0){\line(0,1){90}} \put(20,0){\line(0,1){90}}
\put(40,0){\line(0, 1){90}} \put(80,0){\line(0,1){90}}
\put(100,0){\line(0,1){90}}

\put(7, 100){$E_1$} \put(27, 100){$E_2$} \put(57, 100){$E_3$}
\put(87, 100){$E_4$}

\put (-20, 77){$Q_1$} \put(-20, 50){$Q_2$} \put(-20, 12){$Q_3$}



 \put(90, 63){$\bullet_3$}
 \put(90, 53){$\bullet_4$} \put(90, 43){$\bullet_5$}
 \put(90, 33){$\bullet_6$}  


 \put(1, 63){$\times_3$}
 \put(1, 53){$\times_4$} \put(1, 43){$\times_5$}
 \put(1, 33){$\times_6$} 



\put(200,0){\line(1,0){100}} \put(200,30){\line(1,0){100}}
\put(200,70){\line(1,0){100}} \put(200,90){\line(1,0){100}}

\put(200,0){\line(0,1){90}} \put(220,0){\line(0,1){90}}
\put(240,0){\line(0, 1){90}} \put(280,0){\line(0,1){90}}
\put(300,0){\line(0,1){90}}

\put(207, 100){$E_1$} \put(227, 100){$E_2$} \put(257, 100){$E_3$}
\put(287, 100){$E_4$}

\put(180, 77){$Q_1$} \put(180, 50){$Q_2$} \put(180, 12){$Q_3$}


 \put(290, 23){$\bullet_7$}
 \put(290, 13){$\bullet_8$} \put(290, 3){$\bullet_{9}$}

 \put(221, 23){$\times_7$}
 \put(221, 13){$\times_8$} \put(221, 3){$\times_9$}


\put(21, 90){\line(0,-1){20}} \put(21, 71){\line(1,0){21}}

\put(41.5, 71){\line(0,-1){41}} \put(41, 31){\line(1,0){41}}

\put(81, 31){\line(0,-1){31}}

\put(221, 90){\line(0,-1){20}} \put(221, 71){\line(1,0){21}}

\put(241.5, 71){\line(0,-1){41}} \put(241, 31){\line(1,0){41}}

\put(281, 31){\line(0,-1){31}}


\put(60, -15){$x_{\alpha_p}=\bullet_p$} \put(0,
-15){$x_{\beta_p}=\times_p$}

\put(260, -15){$x_{\alpha_p}=\bullet_p$} \put(200,
-15){$x_{\beta_p}=\times_p$}

\put(140, 50){or} \put(320, 50){etc.}
\end{picture}

\vs \vs \vs

\subsection{Proof of the equality $B_{\ba}=R_{\ba^*}$ }

We will use the same notations as in the previous section.  When
$\{ v_1, \cdots, v_k \}$ is a basis of $\fn_{\ba}$ such that
$v_1=x_{\alpha}$, we will use the notation $x_{\beta} \wedge
\hat{x}_{\alpha} \wedge \cdots \wedge v_k$ to denote the
$k$-vector obtained from $v_1 \wedge \cdots \wedge v_k$, which may
be considered as a base $k$-vector, by replacing $v_1=x_{\alpha}$
with $x_{\beta}$.

 We call the union of all the columns corresponding to $e_p$'s in
$E_i$ the {\it $E_i$-column}. Similarly, we call the union of all
the rows corresponding $e_p$'s in $Q_a$ the {\it $Q_a$-row}.

\begin{prop} \label{equalitygr} Let $\ba=(p_1^{q_1}, \cdots, p_r^{q_r})$
be a partition and let $\ba'=({p'_1}^{q'_1}, \cdots,
{p'_r}^{q'_r})$ be its conjugate. Take a decomposition
$E=\oplus_{i} E_i$ and $Q=\oplus_{a} Q_a$ associated to $\ba$ as
in the previous section. Let $\Pi$ be the index set of $(i,a)$
such that $\fn_{\ba}=\oplus_{(i,a) \in \Pi} E_i^* \otimes Q_a$.
Then the Schubert differential system $\CB_{\ba^*}$ is equal to
the Schur differential system $\CR_{\ba}$ except when both $E_i$
and $Q_a$ are one dimensional for some $(i,a) \not\in \Pi$ with
$E_i^* \otimes Q_a$ adjacent to $\fn_{\ba}$.
\end{prop}

\begin{pf} We will divided the proof by two: I. when the highest
weight of the irreducible component in the complement of
$\fm_{\ba}$ in $\fn_{\ba}^* \otimes \fm/\fn_{\ba}$ is of type (1)
and II. when it is of type (2) or (3) (See Proposition
\ref{typehwv} for the types of the highest weight vectors).

\vs {\bf I.  Type (1)}: Fix $(j,b), \in \Pi, (i,a) \not\in \Pi$.
Let $\{v_1, \cdots, \cdots, v_k \}$ be a basis of $\fn_{\ba}$ with
$v_{1}=x_{\alpha}$ is a lowest weight vector of $E_j^* \otimes
Q_b$ and $x_{\beta}$ is the highest weight vector of $E_i^*
\otimes Q_a$. We will show that $x_{\beta} \wedge \hat{x}_{\alpha}
\wedge \cdots \wedge v_k$ has a nonzero component in $\bI_{\bb}$
for a partition $\bb \not=\ba$ with $|\bb|=|\ba|$.

\vs {\bf Case 1. } Assume that $q_i' \geq 2$ for all $i$. Then we
can take $\bb$ with $\fn_{\bb}=\fn_{\ba} -\{x_{\beta'}\} \cup
\{x_{\alpha'} \}$ where $x_{\beta'}$ is the highest weight in the
boxes $E_j^* \otimes Q_c$ for all $c$ with $(j,c) \in \Pi$ and
$x_{\alpha'}$ is the lowest weight vector in the boxes $E_k^*
\otimes Q_a$ for all $k$ with $(k,a) \not\in \Pi$. Note that there
is $x_{\gamma_j}, j=1,\cdots, 4$ in $\fp_0$ such that $
ad(x_{\gamma_2})ad(x_{\gamma_1}) x_{\beta}=x_{\alpha'} $ and $
ad(x_{\gamma_4})ad(x_{\gamma_3}) x_{\beta'}=x_{\alpha} $ and
$ad(x_{\gamma_j}), j=1, \cdots, 4$

\vs  If $x_{\alpha}, x_{\beta}$, $x_{\beta'}$ and $x_{\alpha'}$
lie neither in the same $E_i$-column nor in the same $Q_a$-raw as
in the picture,

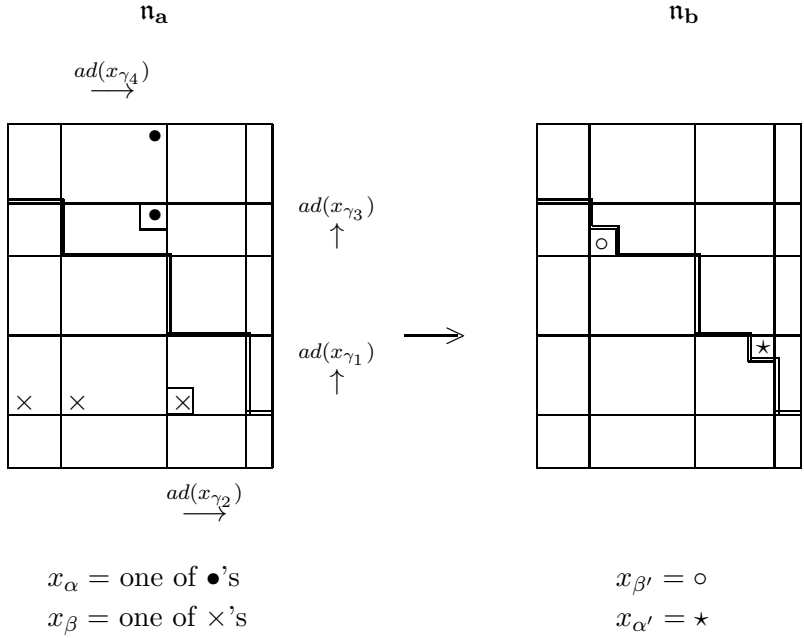
\begin{picture}(300, 200)

\put(50, 170){$\fn_{\ba}$} \put(250, 170){$\fn_{\bb}$}

\put(15, -45){$x_{\alpha}=$ one of $\bullet$'s}

\put(15, -60){$x_{\beta}=$ one of $\times$'s}

\put(230, -45){$x_{\beta'}=\circ$}

\put(230, -60){$x_{\alpha'}=\star$}

\put(25, 140){$\stackrel{ad(x_{\gamma_4})}{\longrightarrow}$}
\put(110, 85){$\stackrel{ad(x_{\gamma_3})}{\uparrow}$}

\put(60, -20){$\stackrel{ad(x_{\gamma_2})}{\longrightarrow}$}
\put(110, 30){$\stackrel{ad(x_{\gamma_1})}{\uparrow}$}

\put(50,90){\usebox{\bx}} \put (60,20){\usebox{\bx}}
\put(53,93){$\bullet$} \put(62, 22){$\times$}
\put(53,123){$\bullet$}

 \put(2,22){$\times$} \put(22,22){$\times$}

\put(0,0){\line(1,0){100}} \put(0,20){\line(1,0){100}}
\put(0,50){\line(1,0){100}} \put(0,80){\line(1,0){100}}
\put(0,100){\line(1,0){100}} \put(0,130){\line(1,0){100}}

\put(0,0){\line(0,1){130}} \put(20,0){\line(0,1){130}}
\put(60,0){\line(0, 1){130}} \put(90,0){\line(0,1){130}}
\put(100,0){\line(0,1){130}}

\put(0,101.5){\line(1,0){21}}
 \put(21, 101.5){\line(0,-1){20}}
\put(20, 81){\line(1,0){41}} \put(61.5,81){\line(0,-1){31}}
\put(60,51){\line(1,0){31}} \put(91.5, 51){\line(0,-1){31}}
\put(90,21.5){\line(1,0){10}}


\put(150,50){\line(1,0){20}} \put(165, 47){$>$}

\put(220,80){\usebox{\bx}} \put (280,40){\usebox{\bx}}
\put(222,82){$\circ$} \put(283, 43){$\star$}

\put(200,0){\line(1,0){100}} \put(200,20){\line(1,0){100}}
\put(200,50){\line(1,0){100}} \put(200,80){\line(1,0){100}}
\put(200,100){\line(1,0){100}} \put(200,130){\line(1,0){100}}

\put(200,0){\line(0,1){130}} \put(220,0){\line(0,1){130}}
\put(260,0){\line(0, 1){130}} \put(290,0){\line(0,1){130}}
\put(300,0){\line(0,1){130}}

\put(200,101.5){\line(1,0){21}} \put(221, 101.5){\line(0,-1){10}}
\put(221, 91.5){\line(1,0){10}} \put(231, 91.5){\line(0,-1){11}}

\put(230, 81){\line(1,0){31}} \put(261.5,81){\line(0,-1){31}}
\put(260,51){\line(1,0){21}} \put(281, 51){\line(0,-1){11}}
\put(281, 41.5){\line(1,0){11}}

\put(291.5, 41){\line(0,-1){21}} \put(290,21.5){\line(1,0){10}}

 \end{picture}

\vs \vs\vs\vs
 \vs\vs\vs

\noindent then  all $\gamma_i$ are distinct. Thus
$$ad(x_{\gamma_4}) \cdots ad(x_{\gamma_1})( x_{\beta} \wedge \cdots
\wedge \hat{x}_{\alpha} \wedge \cdots \wedge v_k)= x_{\alpha'}
\wedge \cdots  \wedge \hat{x}_{\beta'} \wedge \cdots \wedge v_k.
$$ But $x_{\alpha'} \wedge \cdots \wedge
\hat{x}_{\beta'} \wedge \cdots \wedge v_k$ is the lowest weight
vector of $\bI_{\bb}$. Since $ad(\fp_0)$ preserves $\bI_{\bb}$,
there is a nonzero $\bI_{\bb}$-component in $x_{\beta} \wedge
\cdots \wedge \hat{x}_{\alpha} \wedge \cdots \wedge v_k$ and thus
it is not contained in $\bI_{\ba}$.

\vs If  $x_{\alpha}, x_{\beta}$, $x_{\beta'}$ and $x_{\alpha'}$
lie either in the same $E_i$-column or in the same $Q_a$-row as in
the picture

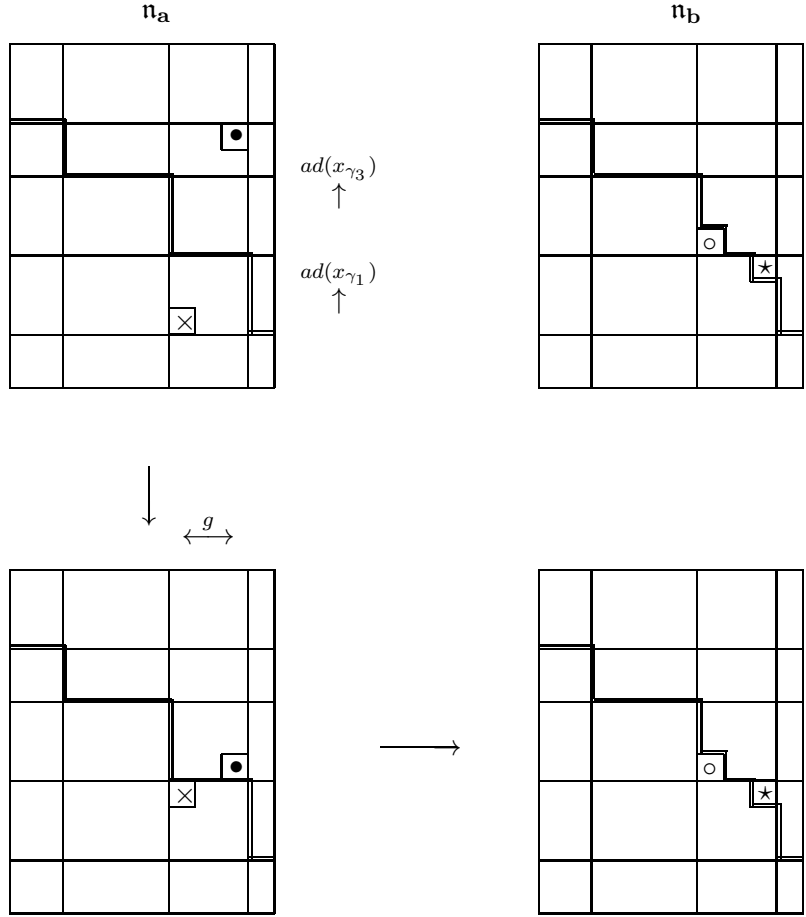
\begin{picture}(300, 170)

\put(50, 140){$\fn_{\ba}$} \put(250, 140){$\fn_{\bb}$}





\put(110, 70){$\stackrel{ad(x_{\gamma_3})}{\uparrow}$}

\put(110, 30){$\stackrel{ad(x_{\gamma_1})}{\uparrow}$}

\put(80,90){\usebox{\bx}} \put (60,20){\usebox{\bx}}
\put(83,93){$\bullet$}

\put(62, 22){$\times$}

\put(0,0){\line(1,0){100}} \put(0,20){\line(1,0){100}}
\put(0,50){\line(1,0){100}} \put(0,80){\line(1,0){100}}
\put(0,100){\line(1,0){100}} \put(0,130){\line(1,0){100}}

\put(0,0){\line(0,1){130}} \put(20,0){\line(0,1){130}}
\put(60,0){\line(0, 1){130}} \put(90,0){\line(0,1){130}}
\put(100,0){\line(0,1){130}}

\put(0,101.5){\line(1,0){21}}
 \put(21, 101.5){\line(0,-1){20}}
\put(20, 81){\line(1,0){41}} \put(61.5,81){\line(0,-1){31}}
\put(60,51){\line(1,0){31}} \put(91.5, 51){\line(0,-1){31}}
\put(90,21.5){\line(1,0){10}}


\put(52.6,-30){\line(0,-1){21}} 
\put(50, -50){$\downarrow$}

\put(260,50){\usebox{\bx}} \put (280,40){\usebox{\bx}}
\put(262,52){$\circ$} \put(283, 43){$\star$}

\put(200,0){\line(1,0){100}} \put(200,20){\line(1,0){100}}
\put(200,50){\line(1,0){100}} \put(200,80){\line(1,0){100}}
\put(200,100){\line(1,0){100}} \put(200,130){\line(1,0){100}}

\put(200,0){\line(0,1){130}} \put(220,0){\line(0,1){130}}
\put(260,0){\line(0, 1){130}} \put(290,0){\line(0,1){130}}
\put(300,0){\line(0,1){130}}

\put(200,101.5){\line(1,0){21}} \put(221, 101.5){\line(0,-1){20}}

\put(221, 81){\line(1,0){41}} \put(261.5,81){\line(0,-1){21}}

\put(261.5, 61.5){\line(1,0){10}} \put(271, 61){\line(0,-1){11}}

 \put(271,51){\line(1,0){11}}
\put(281, 51){\line(0,-1){11}} \put(281, 41.5){\line(1,0){11}}

\put(291.5, 41){\line(0,-1){21}} \put(290,21.5){\line(1,0){10}}

 \end{picture}

\vs \vs\vs\vs

\begin{picture}(300, 150)

\put(65,140){$\stackrel{g}{\longleftrightarrow}$}



\put(80,50){\usebox{\bx}} \put (60,40){\usebox{\bx}}
\put(83,53){$\bullet$}

\put(62, 42){$\times$}

\put(0,0){\line(1,0){100}} \put(0,20){\line(1,0){100}}
\put(0,50){\line(1,0){100}} \put(0,80){\line(1,0){100}}
\put(0,100){\line(1,0){100}} \put(0,130){\line(1,0){100}}

\put(0,0){\line(0,1){130}} \put(20,0){\line(0,1){130}}
\put(60,0){\line(0, 1){130}} \put(90,0){\line(0,1){130}}
\put(100,0){\line(0,1){130}}

\put(0,101.5){\line(1,0){21}}
 \put(21, 101.5){\line(0,-1){20}}
\put(20, 81){\line(1,0){41}} \put(61.5,81){\line(0,-1){31}}
\put(60,51){\line(1,0){31}} \put(91.5, 51){\line(0,-1){31}}
\put(90,21.5){\line(1,0){10}}

\put(140, 62.9){\line(1,0){29}} \put(160, 60){$\rightarrow$}

\put(260,50){\usebox{\bx}} \put (280,40){\usebox{\bx}}
\put(262,52){$\circ$} \put(283, 43){$\star$}

\put(200,0){\line(1,0){100}} \put(200,20){\line(1,0){100}}
\put(200,50){\line(1,0){100}} \put(200,80){\line(1,0){100}}
\put(200,100){\line(1,0){100}} \put(200,130){\line(1,0){100}}

\put(200,0){\line(0,1){130}} \put(220,0){\line(0,1){130}}
\put(260,0){\line(0, 1){130}} \put(290,0){\line(0,1){130}}
\put(300,0){\line(0,1){130}}

\put(200,101.5){\line(1,0){21}} \put(221, 101.5){\line(0,-1){20}}

\put(221, 81){\line(1,0){41}} \put(261.5,81){\line(0,-1){21}}

\put(261.5, 61.5){\line(1,0){10}} \put(271, 61){\line(0,-1){11}}

 \put(271,51){\line(1,0){11}}
\put(281, 51){\line(0,-1){11}} \put(281, 41.5){\line(1,0){11}}

\put(291.5, 41){\line(0,-1){21}} \put(290,21.5){\line(1,0){10}}

 \end{picture}

\vs \vs\vs\vs

\noindent then  $\gamma_2=\gamma_4$ or $\gamma_1=\gamma_3$. We
consider the case when $\gamma_2=\gamma_4$( the proof for the
other case is similar to this case). Thus the multivector
$x_{\alpha'}  \wedge \hat{x}_{\beta'} \wedge \cdots \wedge v_k \in
\bI_{\bb}$ can be obtained from the multivector  $ x_{\beta}
\wedge  \hat{x}_{\alpha} \wedge \cdots \wedge v_k$ by applying $g
\circ ad(x_{\gamma_3})ad(x_{\gamma_1})$, where $g \in SL(m) \times
SL(n-m)$ which exchanges the most left column $e_{s_{i-1}+1}^*
\otimes q_p\,, 1 \leq p \leq n-m$ with the most right column
$e_{s_i}^* \otimes q_p\,, 1 \leq  p \leq n-m$ in the $E_i$-column.
This shows that $ x_{\beta} \wedge \hat{x}_{\alpha} \wedge \cdots
\wedge v_k$ is not contained in $\bI_{\ba}$.

\vs {\bf Case 2. } If some $q_a'=1$, then there may be  no
partition $\bb$ with such property as in Case 1. This is the case
when dim $Q_a$ is one and $ a=b+1$ and $E_{j+1}^* \otimes Q_a
\not\in \fn_{\ba}$. Then $x_{\alpha'}$ is left to $x_{\beta'}$ and
they are adjacent so we cannot find such a partition $\bb$.

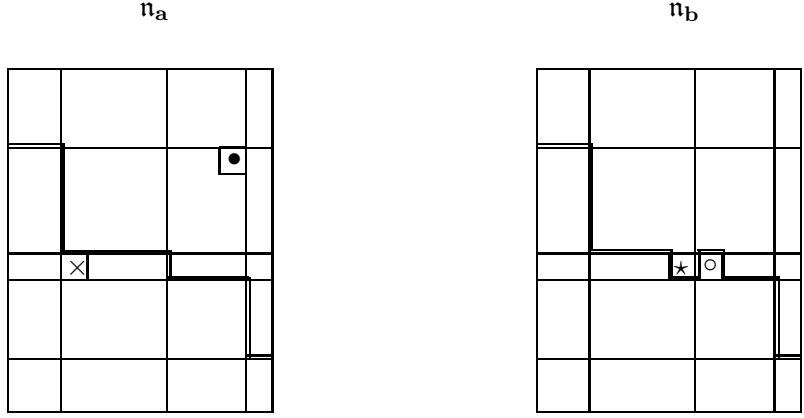
\begin{picture}(300, 180)

\put(50, 150){$\fn_{\ba}$} \put(250, 150){$\fn_{\bb}$}


\put(80,90){\usebox{\bx}} \put (20,50){\usebox{\bx}}
\put(83,93){$\bullet$} \put(22, 52){$\times$}


\put(0,0){\line(1,0){100}} \put(0,20){\line(1,0){100}}
\put(0,50){\line(1,0){100}} \put(0,60){\line(1,0){100}}
\put(0,100){\line(1,0){100}} \put(0,130){\line(1,0){100}}

\put(0,0){\line(0,1){130}} \put(20,0){\line(0,1){130}}
\put(60,0){\line(0, 1){130}} \put(90,0){\line(0,1){130}}
\put(100,0){\line(0,1){130}}

\put(0,101.5){\line(1,0){21}}
 \put(21, 101.5){\line(0,-1){40}}
\put(20, 61){\line(1,0){41}} \put(61.5,61){\line(0,-1){11}}
\put(60,51){\line(1,0){31}} \put(91.5, 51){\line(0,-1){31}}
\put(90,21.5){\line(1,0){10}}

\put(260,50){\usebox{\bx}} \put (250,50){\usebox{\bx}}
\put(263,53){$\circ$} \put(252, 52){$\star$}


\put(200,0){\line(1,0){100}} \put(200,20){\line(1,0){100}}
\put(200,50){\line(1,0){100}} \put(200,60){\line(1,0){100}}
\put(200,100){\line(1,0){100}} \put(200,130){\line(1,0){100}}

\put(200,0){\line(0,1){130}} \put(220,0){\line(0,1){130}}
\put(260,0){\line(0, 1){130}} \put(290,0){\line(0,1){130}}
\put(300,0){\line(0,1){130}}

\put(200,101.5){\line(1,0){21}}
 \put(221, 101.5){\line(0,-1){40}}

\put(220, 61.5){\line(1,0){31}}

\put(251, 61.5){\line(0,-1){11}}

\put(251, 51){\line(1,0){11}}

\put(261.5, 51){\line(0,1){10}}

\put(261, 61.5){\line(1,0){10}}

\put(271, 61){\line(0,-1){11}}

\put(271,51){\line(1,0){21}} \put(291.5, 51){\line(0,-1){31}}
\put(290,21.5){\line(1,0){10}}

\end{picture}
 \vs \vs

But in this case, consider $\bb'$ such that $\fn_{\bb'}=\fn_{\ba}
-\{ x_{\beta''} \} \cup \{ x_{\alpha'} \}$ where $x_{\beta''}$ is
the highest weight vector in the boxes $E_k^* \otimes Q_b$ for all
$k$ with $(k,b) \in \Pi$ and $x_{\alpha'}$ is the lowest weight
vector in the boxes $E_i^* \otimes Q_c$ for all $c$ with $(i,c)
\not\in \Pi$.

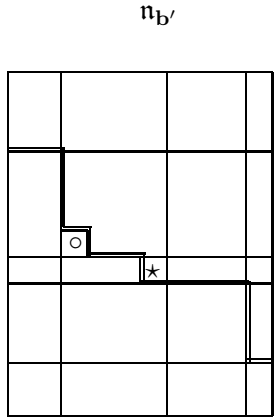
\begin{picture}(300, 180)

\put(50, 150){$\fn_{\bb'}$} 


\put(20,60){\usebox{\bx}} \put (50,50){\usebox{\bx}}
\put(23,63){$\circ$} \put(52, 52){$\star$}


\put(0,0){\line(1,0){100}} \put(0,20){\line(1,0){100}}
\put(0,50){\line(1,0){100}} \put(0,60){\line(1,0){100}}
\put(0,100){\line(1,0){100}} \put(0,130){\line(1,0){100}}

\put(0,0){\line(0,1){130}} \put(20,0){\line(0,1){130}}
\put(60,0){\line(0, 1){130}} \put(90,0){\line(0,1){130}}
\put(100,0){\line(0,1){130}}

\put(0,101.5){\line(1,0){21}}
 \put(21, 101.5){\line(0,-1){30}}

 \put(21, 71.5){\line(1,0){10.5}}

\put(31, 71){\line(0,-1){11}}

\put(31, 61.5){\line(1,0){21}} \put(51.5,61){\line(0,-1){11}}
\put(50,51){\line(1,0){41}} \put(91.5, 51){\line(0,-1){31}}
\put(90,21.5){\line(1,0){10}}
\end{picture}
\vs \vs \vs

 If $i \not=j+1$, then there is such a partition
$\bb'$. If $i=j+1$, then, by the assumption, dim $E_i$ is not
equal to one and thus we can find such a partition $\bb'$. Then by
the same argument as in the case 1, we can prove that $ x_{\beta}
\wedge\hat{x}_{\alpha} \wedge \cdots \wedge v_k$ is not contained
in $\bI_{\ba}$.

\vs {\bf II. Type(2) or (3)} Let $\sum_{q_p \in Q_a}
x_{\alpha_p}^* \otimes x_{\beta_p}$ be the highest weight vector
of type (2) of an irreducible representation space in the
complement of $\fm_{\ba}$ in $\fn_{\ba}^* \otimes \fm/\fn_{\ba}$,
 where $x_{\alpha_p}=e_{r_j }^* \otimes q_p$ and $x_{\beta_p}=e_{r_{i-1}+1 }^*
 \otimes q_p$ for some $(j,a) \in \Pi, (i,a) \not\in \Pi$
 such that $(j, a-1) \in \Pi$ and $(i, a-1) \not \in \Pi$.
We will show that $\sum_{q_p \in Q_a} x_{\beta_p} \wedge
\hat{x}_{\alpha_p} \wedge \cdots \wedge v_k$ has a nonzero
component in $\bI_{\bb} + \bI_{\bc}$  if $\dim Q_a \geq 2$, and
has a nonzero component in $\bI_{\bb}$, otherwise, for some
partition $\bb$ and $\bc$.

Applying the  adjoint actions successively to $$\sum_{q_p \in Q_a}
x_{\beta_p} \wedge \hat{x}_{\alpha_p} \wedge \cdots \wedge
v_k=x_{\beta_{s_{a-1}+1}} \wedge \hat{x}_{\alpha_{s_{a-1}+1} }
\wedge \cdots \wedge v_k
   + \sum_{p =s_{a-1}+2}^{p=s_a} x_{\beta_p} \wedge
\hat{x}_{\alpha_p} \wedge \cdots \wedge v_k,$$

\noindent we  can get
$$ x_{\beta_{s_{a-1}+1}} \wedge  \hat{x}_{\alpha_{s_c} } \wedge
\cdots \wedge v_k \,\,
   +  \,\,x_{\beta_{s_{a-2}+1}} \wedge
\hat{x}_{\alpha_{s_c}} \wedge \cdots \wedge v_k,$$

\noindent where $x_{\beta_{s_{a-2}+1}}=e_{r_{i-1}+1 }^*
 \otimes q_{s_{a-2}+1}$ and $x_{\alpha_{s_c}}=e_{r_j }^* \otimes
 q_{\alpha_{s_c}}$ and $c$ is the largest index $c$ such that $(j,c) \in \Pi$.

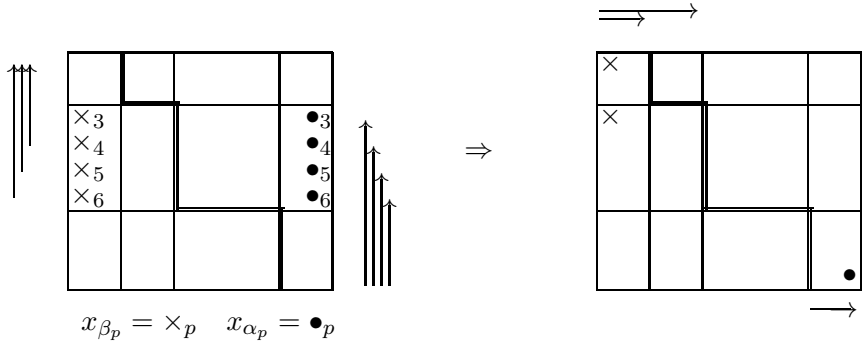
\begin{picture}(300, 130)

\put(0,0){\line(1,0){100}} \put(0,30){\line(1,0){100}}
\put(0,70){\line(1,0){100}} \put(0,90){\line(1,0){100}}

\put(0,0){\line(0,1){90}} \put(20,0){\line(0,1){90}}
\put(40,0){\line(0, 1){90}} \put(80,0){\line(0,1){90}}
\put(100,0){\line(0,1){90}}


\put(110, 57){$\uparrow$} \put(112.5, 2){\line(0,1){60}}

\put(113, 47){$\uparrow$} \put(115.6, 2){\line(0,1){50}}

\put(116, 37){$\uparrow$} \put(118.6, 2){\line(0,1){40}}

\put(119, 27){$\uparrow$} \put(121.6, 2){\line(0,1){30}}

\put(-23, 78){$\uparrow$} \put(-20.4, 35){\line(0,1){50}}

\put(-20, 78){$\uparrow$} \put(-17.4, 45){\line(0,1){40}}

\put(-17, 78){$\uparrow$} \put(-14.4, 55){\line(0,1){30}}


 \put(90, 63){$\bullet_3$}
 \put(90, 53){$\bullet_4$} \put(90, 43){$\bullet_5$}
 \put(90, 33){$\bullet_6$}  


 \put(1, 63){$\times_3$}
 \put(1, 53){$\times_4$} \put(1, 43){$\times_5$}
 \put(1, 33){$\times_6$} 


\put(21, 90){\line(0,-1){20}} \put(21, 71){\line(1,0){21}}

\put(41.5, 71){\line(0,-1){41}} \put(41, 31){\line(1,0){41}}

\put(81, 31){\line(0,-1){31}}


\put(60, -15){$x_{\alpha_p}=\bullet_p$} \put(5,
-15){$x_{\beta_p}=\times_p$}

\put(221, 90){\line(0,-1){20}} \put(221, 71){\line(1,0){21}}

\put(241.5, 71){\line(0,-1){41}} \put(241, 31){\line(1,0){41}}

\put(281, 31){\line(0,-1){31}}





\put(200,0){\line(1,0){100}} \put(200,30){\line(1,0){100}}
\put(200,70){\line(1,0){100}} \put(200,90){\line(1,0){100}}

\put(200,0){\line(0,1){90}} \put(220,0){\line(0,1){90}}
\put(240,0){\line(0, 1){90}} \put(280,0){\line(0,1){90}}
\put(300,0){\line(0,1){90}}


\put(201, 102.7){\line(1,0){15}} \put(208, 100){$\rightarrow$}

\put(201, 105.7){\line(1,0){35}} \put(228, 103){$\rightarrow$}

\put(281, -7.3){\line(1,0){15}} \put(288, -10){$\rightarrow$}

 \put(293, 3){$\bullet$}

 \put(201, 63){$\times$}
 \put(201, 83){$\times$}


\put(21, 90){\line(0,-1){20}} \put(21, 71){\line(1,0){21}}

\put(41.5, 71){\line(0,-1){41}} \put(41, 31){\line(1,0){41}}

\put(81, 31){\line(0,-1){31}}

\put(221, 90){\line(0,-1){20}} \put(221, 71){\line(1,0){21}}

\put(241.5, 71){\line(0,-1){41}} \put(241, 31){\line(1,0){41}}

\put(281, 31){\line(0,-1){31}}


\put(150, 50){$\Rightarrow$}


\end{picture}

\vs \vs \vs\vs



Applying the adjoint actions again, we can get

$$ y_{\beta_{s_{a-1}+1}} \wedge  \hat{y}_{\alpha_{s_c} } \wedge
\cdots \wedge v_k \,\,
   +  \,\,y_{\beta_{s_{a-2}+1}} \wedge
\hat{y}_{\alpha_{s_c}} \wedge \cdots \wedge v_k,$$

\noindent where $y_{\beta_{s_{a-1}+1}}=e_{r_{h}}^*
 \otimes q_{s_{a-1}+1}$ and
$y_{\beta_{s_{a-2}+1}}=e_{r_{i}}^*
 \otimes q_{s_{a-2}+1}$ and $y_{\alpha_{s_c}}=e_{r_{j-1}+1 }^* \otimes
 q_{\alpha_{s_c}}$ and $h$ is the largest index $h$ such that $(h,a) \not \in\Pi$.
 This is the sum of the lowest weight vector of $\bI_{\bb}$ and
that of $\bI_{\bc}$ for some partition $\bb$ and $\bc$ such that
all $\ba$, $\bb$, $\bc$ are distinct.

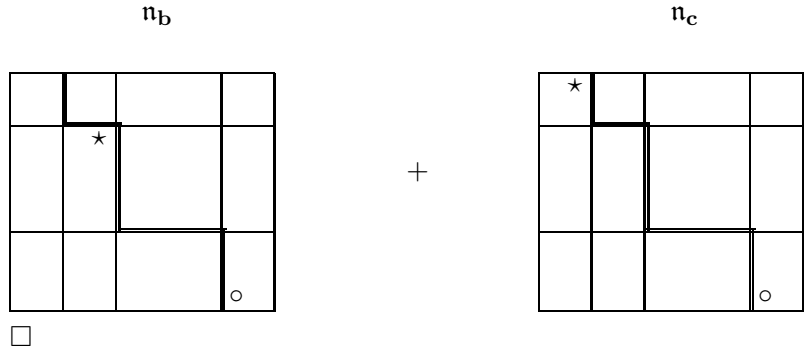
\begin{picture}(300, 130)


\put(50, 110){$\fn_{\bb}$} \put(250, 110){$\fn_{\bc}$}


\put(0,0){\line(1,0){100}} \put(0,30){\line(1,0){100}}
\put(0,70){\line(1,0){100}} \put(0,90){\line(1,0){100}}

\put(0,0){\line(0,1){90}} \put(20,0){\line(0,1){90}}
\put(40,0){\line(0, 1){90}} \put(80,0){\line(0,1){90}}
\put(100,0){\line(0,1){90}}



 \put(83, 3){$\circ$}

 \put(31, 63){$\star$}



\put(200,0){\line(1,0){100}} \put(200,30){\line(1,0){100}}
\put(200,70){\line(1,0){100}} \put(200,90){\line(1,0){100}}

\put(200,0){\line(0,1){90}} \put(220,0){\line(0,1){90}}
\put(240,0){\line(0, 1){90}} \put(280,0){\line(0,1){90}}
\put(300,0){\line(0,1){90}}



 \put(211, 83){$\star$}
\put(283, 3){$\circ$}

\put(21, 90){\line(0,-1){20}} \put(21, 71){\line(1,0){21}}

\put(41.5, 71){\line(0,-1){41}} \put(41, 31){\line(1,0){41}}

\put(81, 31){\line(0,-1){31}}

\put(221, 90){\line(0,-1){20}} \put(221, 71){\line(1,0){21}}

\put(241.5, 71){\line(0,-1){41}} \put(241, 31){\line(1,0){41}}

\put(281, 31){\line(0,-1){31}}



\put(150, 50){$+$} 
\end{picture}

\end{pf}

By Theorem \ref{Schubertgr} and  Proposition \ref{equalitygr}, we
get

\begin{thm} \label{Schurgr}
Let $\ba=(p_1^{q_1}, \cdots, p_r^{q_r}), p_r \not=0$ be a
partition and let $\ba'=({p'_1}^{q'_1}, \cdots, {p'_r}^{q'_r})$, $
p'_r \not=0$ be its conjugate.  Then $\sigma_{\ba}$ is Schur rigid
if
 $q_i, q'_i \geq 2$ for all $i \leq r$.
\end{thm}

\begin{rk}

One of the problems in algebraic geometry is the smoothability of
a singular Schubert variety $X_w$ of $G/P$. We say $X_w$ is {\it
smoothable} if  there is a smooth subvariety $X$ of $G/P$ with
$[X]=[X_w]$ in $H_*(G/P, \BZ)$(\cite{B}). Assume that $\ba=(p^q)$
and that $p=1$ or $q=1$ but both are not 1. Then $\sigma_{\ba}$ is
a singular Schubert variety and the non-smoothability of $X_{\ba}$
is proved in \cite{B}: if $X$ is a subvariety of $Gr(m,n)$ with
$[X]=[\sigma_{\ba}]$, then $X$ is a Schubert variety of type
$\ba$. By Theorem \ref{Schurgr}, for a partition $\ba=(p_1^{q_1},
\cdots, p_r^{q_r})$ with its conjugate $\ba'=({p'_1}^{q'_1},
\cdots, {p'_r}^{q'_r})$, if $q_i, q'_i \geq 2$, for all $i$, then
the singular Schubert varieties $\sigma_{\ba}$ of type $\ba$ is
not smoothable, neither.
\end{rk}

\begin{rk} With the same notations as in Proposition
\ref{equalitygr}, if both $E_i$ and $Q_a$ are one dimensional for
some $(i,a) \not\in \Pi$ with $E_i^* \otimes Q_a$ adjacent to
$\fn_{\ba}$, then $B_{\ba}$ is a proper subvariety of $R_{\ba^*}$.

Consider the highest weight vector $x_{\alpha}^* \otimes
x_{\beta}$ of an irreducible component in the complement of
$\fm_{\ba}$ in $\fn_{\ba}^* \otimes \fm/\fn_{\ba}$ such that
$x_{\alpha} \in \fn_{\ba}$ and $x_{\beta} \in E_i^* \otimes Q_a
\subset \fm/\fn_{\ba}$. Then one can check that $x_{\beta} \wedge
\hat{x}_{\alpha} \wedge \cdots \wedge v_k$ is contained in
$\bI_{\ba}$. Thus this gives a nontrivial element in $
T_{\ba}:=T_{[\fn_{\ba}]}Gr(k, \fm) \cap T_{\wedge^k \fn_{\ba}}
\PP(\bI_{\ba}) \subset \fn_{\ba}^* \otimes \fm/\fn_{\ba}$. Note
that as an element of $T_{[\wedge^k \fn_{\ba}] }\PP(\bI_{\ba})$,
this tangent vector gives the map
$$x_{\alpha} \wedge v_2 \wedge \cdots \wedge v_k \longmapsto
x_{\beta} \wedge v_2 \wedge \cdots \wedge v_k \in \bI_{\ba} /
\wedge^k \fn_{\ba},
$$ \noindent where $\{ v_1=x_{\alpha}, v_2, \cdots, v_k \}$ be a
basis of $\fn_{\ba}$. Then $c(t)=(x_{\alpha} + t x_{\beta}) \wedge
v_2 \wedge \cdots \wedge v_k$, $t \in \BC$, is a curve in
$R_{\ba}=Gr(k, \fm) \cap \PP(\bI_{\ba})$ whose tangent vector is
$x_{\alpha}^* \otimes x_{\beta}$. But this tangent vector is not
contained in $\fm_{\ba}$ and thus $B_{\ba}$ is a proper subvariety
of $R_{\ba^*}$.
 This is a generalization of the counterexamples considered in
\cite{W} or Example 9 of \cite{B}.
\end{rk}

 \vs

\vs\vs

Research Institute of Mathematics

Seoul National University

San 56-1 Sinrim-dong Kwanak-gu

Seoul, 151-747 Korea

jhhong@math.snu.ac.kr



\begin{thebibliography}{99}


\bibitem[B]{B} R. Bryant, \emph{Rigidity and quasi-rigidity of extremal
cycles in compact Hermitian symmetric spaces}, Annals of
mathematics studies 153, Princeton University Press, 2004



\bibitem[FH]{FH} W. Fulton and J. Harris, Representation Theory; A
First Course, Springer-Verlag, 1991




\bibitem[G]{G} A. B. Goncharov, \emph{Generalized conformal structures
on manifolds}, Selecta Mathematica Sovietica, vol 6, No 4, 1987.


\bibitem[H]{H} J. Harris, Algebraic geometry, GTM 133, Springer-Verlag,
1992

\bibitem[Ho]{Ho} J. Hong, Rigidity of smooth Schubert varieties in
Hermitian symmetric spaces, math.DG/0410138






\bibitem[W]{W} M. Walters, \emph{Geometry and uniqueness of some
extreme subvarieties in complex Grassmannians}, Ph.D. thesis,
University of Michigan, 1997.


\end{thebibliography}
\end{document}